\documentclass[11pt]{article}

\usepackage{amsmath,amsthm,amsfonts,amssymb,amscd}
\usepackage{mathtools}
\usepackage{pgf,tikz}
\usepackage{float}
\usepackage{enumerate}
\usepackage{booktabs}
\usepackage{pgfplots}
\usepackage[noend]{algpseudocode}
\usepackage{graphicx,bm,xcolor}
\usepackage{algorithm}
\usepackage{algpseudocode}
\usepackage{hyperref}
\usepackage{caption}
\usepackage[T1]{fontenc}
\usepackage{t1enc}
\usepackage[polish,german,english]{babel}
\usepackage{verbatim} 

\usepackage{url}
\usepackage{authblk}

\usepackage[capitalize]{cleveref} % TW: for cref. MCS: capitalize.
\usepackage{siunitx}
\usepackage{booktabs}
\usepackage{pgfplotstable}
\usepackage{stmaryrd} % Where \llbracket and \rrbracket are defined
\usepackage[top=2cm, bottom=2cm, left=2cm, right=2cm]{geometry}
\newtheorem{problem}{Problem}
\theoremstyle{definition}

\newtheorem{remark}{Remark}

\normalsize

\newcommand\iprod[3][]{\bigl(#2,#3\bigr)_{#1}}% Inner product (#2,#3)_#1
\newcommand\prodf{\iprod[\Omega_f]}
\newcommand\prods{\iprod[\Omega_s]}
\newcommand\prodi{\iprod[\Gamma_i]}
\newcommand\prodd{\iprod[\Gamma_f \setminus \Gamma_f^D]}
\newcommand\norm[1]{\lVert#1\rVert}

\newcommand\jump[1]{\llbracket#1\rrbracket}% [[#1]]

\newcommand\jumpn[2]{\jump{\partial_n^{#2}\varphi_{#1}}}
\newcommand\Jump[1]{g_F^{#1}(\varphi_1, \varphi_2)}
\newcommand\fcdot{\,\cdot\,}
\newcommand\bd{\partial}% Boundary
\newcommand\cl{\overline}% Closure: maybe replace with \bar
\newcommand\tsb[1]{_{\textup{#1}}}% Text subscript in math mode
\newcommand\set[1]{\{#1\}}
\newcommand\defset[2]{\set{#1 \ | \ #2}}
\newcommand\R{\mathbb R}
\renewcommand\P{\mathcal P}% Pilcrow not needed.
\newcommand\T{\mathcal T}
\newcommand\U{\mathcal U}
\newcommand\V{\mathcal V}
\newcommand\X{\mathcal X}
\newcommand\wmax{w\tsb{max}}
\newcommand\tria[1]{\T_{h,n}^{#1}}
\newcommand\triaext[1]{\T_{h,n}^{#1,\textup{ext}}}
\newcommand\Oext[1]{\Omega_{h,n}^{#1,\textup{ext}}}
\newcommand\Ocomp[1]{\Omega_{h,n}^{#1}}
\newcommand\FG[1]{\mathcal{F}_G^{#1}}
\newcommand\FExt[1]{\mathcal{F}_{h,n}^{#1,\textup{ext}}}
\let\0\emptyset% oder \varnothing: das ging durcheinander.
\DeclareMathOperator\tr{tr}
\newcommand\ghw[3][]{g^{h,n,w}\ifx#1\empty\else_{#1}\fi\left(#2,#3\right)}
\newcommand\ghwext[3][]{\ifx#1\empty\ghw{#2}{#3}\else\ghw[#1,\textup{ext}]{#2}{#3}\fi}
\newcommand\gammaext[1]{\ensuremath{\gamma^{\textup{ext}}_{#1}}}

\newcommand\intI{I_T}

\newcommand\restr[2]{{% we make the whole thing an ordinary symbol
  \left.\kern-\nulldelimiterspace % automatically resize the bar with \right
  #1 % the function
  \vphantom{\big|} % pretend it's a little taller at normal size

  \right|_{#2} % this is the delimiter
  }}

\begin{document}

\title{Numerical Simulations of Fully Eulerian Fluid-Structure Contact Interaction using a Ghost-Penalty Cut Finite Element Approach}
\author[1]{S. Frei}
\author[2]{T. Knoke}
\author[2]{M.C. Steinbach}
\author[2]{A.-K. Wenske}
\author[2]{T. Wick}

\affil[1]{University of Konstanz, Department of Mathematics \& Statistics, Universit\"atsstraße 10, 78457 Konstanz, Germany}
\affil[2]{Leibniz Universit\"at Hannover, Institut f\"ur Angewandte
  Mathematik, Welfengarten 1, 30167 Hannover, Germany}

\date{}

\maketitle
	
%%%%%%%%%%%%%%%%%%%%%%%%%%%%%%%%%%%%%%%%%%%%%%%%%%%%%%%%%%%%	
\begin{abstract}
In this work, we develop a
cut-based unfitted finite element formulation
for solving nonlinear, nonstationary fluid-structure interaction with contact in Eulerian coordinates.
In the Eulerian description fluid flow modeled by the incompressible
Navier-Stokes equations remains in Eulerian coordinates, while elastic solids
are transformed from Lagrangian coordinates into the Eulerian system.
A monolithic description is adopted.
For the spatial discretization,
we employ an unfitted finite element method with ghost penalties
based on inf-sup stable finite elements. To handle contact, we use a
relaxation of the contact condition in combination with a unified Nitsche
approach that takes care implicitly of the switch between fluid-structure
interaction and contact conditions. The temporal discretization is based on a
backward Euler scheme with implicit extensions of solutions at the previous time step.
The nonlinear system is solved with
a semi-smooth Newton's method with line search.
Our formulation, discretization and implementation
are substantiated with an elastic falling ball that comes into
contact with the bottom boundary, constituting a challenging state-of-the-art
benchmark.
\end{abstract}

\newpage
%%%%%%%%%%%%%%%%%%%%%%%%%%%%%%%%%%%%%%%%%%%%%%%%%%%%%%%%%%%%
\section{Introduction}
\label{sec_intro}
Fluid-structure interaction is a prime example for a multiphysics
problem. The main challenges are concerned with the interface coupling
and the different coordinate systems of fluids and solids.
Numerous works with different methodologies have been published,
such as the arbitrary Lagrangian Eulerian method~\cite{DoneaSurvey, Hughes1981, FoNo99}, fictious
domain methods~\cite{Glowinski1994283,Glowinski2001363}, immersed boundary methods~\cite{Pe02},
and fully Eulerian methods~\cite{Du06,CoMaMi08}. Subsets and further
refinements of these techniques exist. Classical monographs
and textbooks on fluid-structure interaction
include~\cite{FoQuaVe09,BoGaNe14,GaRa10,richter2017book,FrHoRiWiYa17,BaTaTe13}.

In this paper, we concentrate on fully Eulerian fluid-structure
interaction. As mentioned above, it has been introduced
in~\cite{Du06,CoMaMi08} and has since then been investigated and improved in several studies, such as \cite{richterWick2010,SuIiTaTaMa11,Wi12_fsi_euler,RICHTER2013227,LaRuiQua13,
Frei2015EulerianTechniquesI,Sun20141,RATH2023112188}.
The benefit of a fully Eulerian formulation for FSI
lies in its ability to handle naturally very large deformations,
topology changes, and contact problems in a straightforward way,
see e.g.~\cite{FrRiWi16_JCP,FrRi17,HechtPironneau2016,BurmanFernandezFrei2020}.
The idea is to formulate both the flow
and the solid problem in Eulerian coordinates in time-dependent domains $\Omega_f(t)$ resp.~$\Omega_s(t)$. An accurate numerical method requires the resolution of the interface $\Gamma_i(t)$ separating $\Omega_f(t)$ and $\Omega_s(t)$, which can move freely depending on the solid displacements.
The construction of a fitted finite element method is cumbersome when the interface moves, see e.g.~\cite{FreiRichter2014, FreiPressure2019}. An elegant alternative is given by the cut finite element method~\cite{HANSBO2002elliptic, hansbo2005Nitsche, Burman2015CutFEM} which is based on a fixed finite
element mesh for all times.

In cut finite element methods (CutFEM) \cite{Burman2015CutFEM}, interface conditions are imposed by means of Nitsche's method \cite{Nitsche1971},
see also \cite{HANSBO2002elliptic, hansbo2005Nitsche}. Moreover, additional stabilization via ghost penalty terms at the faces of cut cells is proposed, since the condition number of the system matrix suffers from cells cut into vastly
different sizes, see also \cite{BURMAN2010GhostPenalty, BURMAN2012FictitiousDomainII}. This adaptation of Nitsche's method is used in \cite{BURMAN2014FSI}, where
linear Stokes flow is coupled to a linear elastic structure through separate overlapping meshes, where the solid is described in Lagrangian coordinates on a fitted mesh and glued to the (unfitted) fluid mesh. Fictitious domain methods using cut elements with stabilized Nitsche's method for Stokes’ problem
were investigated in \cite{burman_hansbo_2014}.

This work is an extension of the proceedings paper \cite{FrKnStWeWi24_ENUMATH}, where we employ a cut finite element method for implementing a variational-monolithic fully Eulerian fluid-structure interaction formulation
on a fixed single mesh and without contact. While \cite{FrKnStWeWi24_ENUMATH}
was on a quasi-stationary setting, we perform in the current work
conceptional developments towards nonstationary settings
with moving interfaces.
Specifically, the solid convection terms are implemented and stabililized
with the help of SUPG (streamline upwind Petrov-Galerkin)~\cite{BrHu82}.
Additional challenges arise in the time discretization of the equations that are posed on moving subdomains $\Omega_f(t), \Omega_s(t)$. As $\Omega_l(t_n) \neq \Omega_l(t_{n-1})$ ($l\in \{f,s\}$), solutions from the previous time step need to be extended before they can be used within a time-stepping scheme on $\Omega_l(t_n)$. Here, we follow the approach introduced by Lehrenfeld \& Olshanskii in~\cite{LehOl19} for parabolic equations on time-dependent domains, where the solutions are already extended implicitly in the previous time step by means of ghost-penalty extensions. While the approach has been introduced for BDF-type time-stepping schemes in~\cite{LehOl19}, a Crank-Nicolson-type scheme has recently been analysed in~\cite{FreiSingh2024}. Corresponding BDF schemes for the Stokes equations have been introduced and analysed in~\cite{burman2022eulerian,vonunfitted21}.
The resulting formulation is treated
in a monolithic fashion by using a line-search Newton method.
Therein, the linear equation
systems are solved by means of MUMPS~\cite{MUMPS:1},
i.e., a parallel sparse direct solver.

Recently, the combination of fluid-structure interactions and contact has been the subject of several research papers~\cite{TezduyarSathe2007, AstorinoGerbeauetal2009, AgerWalletal, BurmanFernandezFrei2020, vonWahletal2021, BurmanFernandezFreiGerosa2024, FormaggiaGattiZonca2021, FaraSchwarzacherTuma2024, GerosaMarsden2024}. When the solid comes into contact with another solid or an exterior wall, several numerical challenges need to be tackled, including a topology change in the fluid domain, the Navier-Stokes no-collision paradox~\cite{Hillairet2d, HeslaPhD, HillairetTakahashi3d} and a non-smooth switch from fluid-structure interface to contact conditions. A simple way to circumvent these difficulties is a relaxation of the contact conditions~\cite{BurmanFernandezFrei2020, BurmanFernandezFreiGerosa2024, FormaggiaGattiZonca2021}, where the no-penetration condition is already imposed at a distance of $\epsilon>0$ between the solids, such that a small fluid layer of size $\epsilon>0$ remains for all times. Different extensions have been proposed, in order to give a physical meaning to the fluid layer, considering rough surfaces of the solids~\cite{AgerWalletal, BurmanFernandezFreiGerosa2024, GerosaMarsden2024}. If a fluid layer is present for all times, the method of Alart and Curnier~\cite{AlartCurnier1991} can be used to reformulate contact and fluid-structure interaction conditions into a system of (non-smooth) equalities and both conditions can be incorporated in an implicit way by means of a unified Nitsche approach~\cite{BurmanFernandezFrei2020}.

Finally, our newly developed approach and its algorithmic
realization are implemented in the finite element software
deal.II~\cite{dealii2019design,dealII95}. Therein,
as numerical example
an elastic falling ball is considered,
which is part of benchmark computations~\cite{vWRFH21}, and
a challenging application. As quantities of interest,
we observe
distances from the bottom, contact times, impact velocities, elastic and kinetic energies, and performance
of the nonlinear Newton solver.

The outline of this paper is as follows. In \Cref{sec_FSI} our
fully Eulerian fluid-structure interaction formulation is presented.
Then, in~\Cref{sec_disc} we introduce the discretization in space and time, including implicit extensions. In~\Cref{sec_cont}, we introduce the FSI-contact problem and our numerical approach to handle contact. Next, we derive a semi-smooth Newton method in~\Cref{sec_sol}. Finally, our algorithm is substantiated by numerical examples in~\Cref{sec_tests}.

\newpage
%%%%%%%%%%%%%%%%%%%%%%%%%%%%%%%%%%%%%%%%%%%%%%%%%%%%%%%
\section{Fully Eulerian Fluid-Structure Interaction}
\label{sec_FSI}
In this section, we introduce a fully Eulerian fluid-structure interaction
formulation. First, the
strong formulation is introduced, followed by its weak counterpart.

\subsection{Strong Form}
Let $\Omega \subseteq \R^2$ be a fixed bounded domain, which is partitioned into a
fluid subdomain $\Omega_f:=\Omega_f(t)$ and a
solid subdomain $\Omega_s:=\Omega_s(t)$ such that
$\cl\Omega = \cl\Omega_f(t) \cup \cl\Omega_s(t)$ with
$\Omega_f(t) \cap \Omega_s(t) = \0$ for all $t \in \intI:=[0,T_\mathrm{end}]$.
In an Eulerian description, the solid sub-domain $\Omega_s(t)$ is defined implicitly by the (unknown) solid displacement $u$ as
\begin{align}\label{MovingDomains}
 \Omega_s(t) = \big\{ x\in\Omega  \, \big| \, T(x,t) \in \Omega_s(0)\big\},
\end{align}
where $T\colon\Omega(t)\to\Omega$ is a bijective mapping that is given by $T(x,t) =x-u(x,t)$
in the solid domain $\Omega_s(t)$ and by an arbitrary (smooth) extension in $\Omega\setminus \Omega_s(t)$.
For the details, we refer to the textbook~\cite{richter2017book}. The fluid sub-domain $\Omega_f(t)$ is then defined by the interior of $\Omega\setminus \Omega_s(t)$, i.e. $\Omega_f(t):= \text{int}(\Omega\setminus \Omega_s(t))$.
We assume that both $\Omega_f(t)$ and $\Omega_s(t)$
are parameterized by a $C^{0,1}$ boundary,
such that all terms arising in the following equations
are well-defined at any time.

Next, let $\Gamma_i := \Gamma_i(t) := \cl\Omega_f(t) \cap \cl\Omega_s(t)$
be the interface between the subdomains, and
$\Gamma_f^D := \Gamma_f^D(t) \subseteq \Gamma_f(t) := \bd\Omega_f(t) \cap \bd\Omega$ and
$\Gamma_s^D := \Gamma_s^D(t) \subseteq \Gamma_s(t) := \bd\Omega_s(t) \cap \bd\Omega$ the Dirichlet boundaries. Finally, we write
\begin{align*}
Q_l&:=\bigcup_{t\in [0,T_\mathrm{end}]} \Omega_l(t) \cup \{t\},& \qquad \Sigma_l^D&:=\bigcup_{t\in [0,T_\mathrm{end}]} \Gamma_l^D(t) \cup \{t\},
\quad l\in\{s,f\},\\
Q&:=\Omega \times I_T,\qquad& \Sigma_l&:=\bigcup_{t\in [0,T_\mathrm{end}]} \Gamma_l(t) \cup \{t\},
\quad l\in\{s,f,i\},
\end{align*}
for the space-time subdomains.

For the fluid velocity $v_f\colon \cl Q_f \to \R^2$,
the pressure $p\colon \cl Q_f \to \R$,
the solid velocity $v_s\colon \cl Q_s \to \R^2$ and
the displacement $u\colon \cl Q_s \to \R^2$,
we define stresses
\begin{align}
 \sigma_f &:= \sigma_f(v_f,p) := \rho_f \nu_f (\nabla v_f + \nabla v_f^\top) - pI \label{fstress} \\
  \text{and}\quad \sigma_s &:= \sigma_s(u) := 2 \mu_s E_s + \lambda_s \tr(E_s)I,
\end{align}
where $E_s := \frac12 (\nabla u + \nabla u^\top + \nabla u^\top \cdot \nabla u)$
denotes the nonlinear Green-Lagrange strain,
$\mu_s$ and $\lambda_s$ are the Lam\'{e} parameters,
$\rho_f$ and $\rho_s$ are the densities of the fluid and the solid,
and $\nu_f$ is the fluid viscosity. For the solid , this corresponds to the
nonlinear St.\ Venant-Kirchhoff material law.

Moreover, $f\colon Q \to \R^2$
is a given external force field,
$v_f^D\colon \Sigma_f^D \to \R^2$
and
$u^D\colon \Sigma_s^D  \to \R^2$
are functions on the Dirichlet boundaries, and
$\smash{v_f^0}\colon \Omega_f(0) \to \R^2$,
$v_s^0\colon \Omega_s(0) \to \R^2$ and
$u^0\colon \Omega_s(0) \to \R^2$ describe initial values.

The fully Eulerian FSI system is then defined as follows in strong form.
\begin{problem}
\label{problem_strong}
Find $(v_f,p,v_s,u)$ such that
\begin{align*}
&\left\{\begin{array}{lll}
        \rho_f \partial_t v_f+\rho_f (v_f \cdot \nabla)v_f-\nabla \cdot \sigma_f &=\rho_f f \quad & \text{in } Q_f, \\
        \nabla \cdot v_f &=0 \quad & \text{in } Q_f, \\
        v_f &=v_f^D \quad & \text{on } \Sigma_f^D, \\
        \rho_f \nu_f \partial_n v_f-p n &= 0%- \rho_f \nu_f \nabla v_f^\top\cdot n \quad
        & \text{on } \Sigma_f \setminus \Sigma_f^D, \\
        v_f &= v_f^0 & \text{in } \Omega_f(0),
       \end{array} \right. \\[-0.5pt]
&\left\{\begin{array}{lll}
        \rho_s \partial_t v_s +\rho_s (v_s \cdot \nabla)v_s-\nabla \cdot \sigma_s &=\rho_s f \quad & \text{in } Q_s, \\
        \partial_t u + (v_s \cdot \nabla)u - v_s &=0 \quad & \text{in } Q_s, \\
        u &=u^D \quad & \text{on } \Sigma_s^D, \\
        \sigma_s \cdot n &=0 \quad & \text{on } \Sigma_s  \setminus \Sigma_s^D, \\
        u &= u^0 & \text{in } \Omega_s(0), \\
        v_s &= v_s^0 & \text{in } \Omega_s(0),
       \end{array} \right. \\[-0.5pt]
&\left\{\begin{array}{lll}
        v_f &= v_s \quad & \text{on } \Sigma_i, \\
        \sigma_f \cdot n &=\sigma_s \cdot n \quad & \text{on }
\Sigma_i.
       \end{array} \right.
\end{align*}
\end{problem}

%%%%%%%%%%%%%%%%%%%%%%%%%%%%%%%%%%%%%%%%%%%%%%%%%%%%%%%%%%%%%
\subsection{Weak Formulation}

Let $\V_f := H_0^1(\Omega_f;\Gamma_f^D)$,
$\V_s := H^1(\Omega_s)$, 
$\U := H_0^1(\Omega_s;\Gamma_s^D)$ and
$\P := L^2(\Omega_f)$ be given function spaces.
Due to the outflow condition on
$\Gamma_f \setminus \Gamma_f^D$, this definition for the pressure space $\P$ is sufficient to guarantee uniqueness of the pressure.
We denote the product space by
$$\X := \V_f \times \P \times \V_s \times \U. $$
For a function space $V$, let
$W(\intI, V, V^*) := \defset{v \in L^2(\intI,V)}{\partial_t v \in L^2(\intI,V^*)}$,
where $L^2(\intI,\fcdot)$ are Bochner spaces and $V^*$ denotes the dual space of $V$.
We combine the Bochner spaces into
$$W(\intI, \X) := W(\intI,\V_f,\V_f^*) \times L^2(\intI,\P) \times W(\intI, \V_s, \V_s^*) \times W(\intI, \U, L^2(\Omega)).$$
Now, based on Problem~\ref{problem_strong}, we can state the weak formulation of the FSI system as follows.
\begin{problem}
Find
$(v_f, p, v_s, u) \in (v_f^D,0,0,u^D) + W(\intI,\X)$
such that
for almost all $t \in \intI$ it holds
$v_f = v_s$ on $\Gamma_i(t)$ and for all $(\phi_f, \psi, \phi_s, \xi) \in \X$:
\begin{multline*}
  \rho_f \langle\partial_t v_f + v_f \cdot \nabla v_f,\phi_f\rangle_{V_f^*\times V_f} + \prodf{\sigma_f}{\nabla\phi_f} + \prodf{\nabla\cdot v_f}{\xi}
  \\
  - \langle{\rho_f \nu_f \nabla v_f^\top n_f, \phi_f\rangle_{H^{-1/2}(\Gamma_i)\times H^{1/2}(\Gamma_i)}} - \langle\sigma_f\cdot n_f,\phi_f - \phi_s\rangle_{H^{-1/2}(\Gamma_i)\times H^{1/2}(\Gamma_i)}\\
  + \rho_s \langle \partial_t v_s + v_s\cdot \nabla v_s, \phi_s\rangle_{V_s^*\times V_s} + \prods{\sigma_s}{\nabla\phi_s}
    +\langle\partial_t u + v_s \cdot \nabla u - v_s,\psi\rangle_{V_s^*\times V_s}\\
    = \rho_f \prodf{f}{\phi_f} + \rho_s \prods{f}{\phi_s},
  \end{multline*}
  where $\langle\fcdot,\fcdot\rangle_{V\times V^*}$ denote respective duality pairings of $V$ and its dual space $V^*$.
\end{problem}

%%%%%%%%%%%%%%%%%%%%%%%%%%%%%%%%%%%%%%%%%%%%%%%%%%%%%%%
\section{Discretization using Ghost Penalities and Cut Finite Elements}
\label{sec_disc}
In this section, discretizations in time and space by a cut Finite Element (FE)
approach are introduced.

\subsection{Temporal Discretization}
We follow the Rothe method and perform first discretization in time
and then discretization in space.
For the temporal discretization, we utilize the A-stable
first-order backward Euler method
with time grid $0 = t_0 < \dots < t_N = T_\mathrm{end}$.
For simplicity, we will consider equidistant time steps $k=t_{n+1}-t_n$. As the interface moves with time, we will need to extend the solution at time $t_n$ to a larger domain that contains both $\Omega_l(t_n)$ and $\Omega_l(t_{n+1})$ ($l=s,f$), where it will be needed in the following time step.
Details on the extension will be given in Subsection~\ref{subsec.ext}.

\subsection{Spatial Discretization with Ghost Penalities and Cut Finite Elements}
In the spatial discretization we use
continuous quadratic elements for the fluid velocity and
continuous linear elements for the remaining solution components.
Let $\T_h$ be a quasi-uniform triangulation of $\Omega$
that is fitted to the boundary of the domain $\Omega$,
but not to the interface $\Gamma_i$.
Moreover, let
\begin{align*}
  \tria{f} := \defset{T \in \T_h}{T \cap \Omega_f(t_n) \ne \0}
  \quad\text{and}\quad
  \tria{s} := \defset{T \in \T_h}{T \cap \Omega_s(t_n) \ne \0}
\end{align*}
be overlapping sub-triangulations.
We use the following finite element spaces on $\tria{l}$:
\newcommand\Vhn[2]{V_{h,n,#1}^{(#2)}}
\begin{align*}
  \Vhn{l}{r} :=
  \defset{\phi \in C(\cl{\Ocomp{l}})}
  {\phi|_T \in Q_r(T) \ \forall T \in \tria{l}},
  \quad l \in \set{f, s},
\end{align*}
where $\Ocomp{l}$ denotes the domain spanned by the cells $T\in \tria{l}$.
We define spaces $\V_{h,{n,}f}:=\Vhn{f}{2} \cap H_0^1(\Omega_f;\Gamma_f^D)$,
$\V_{h,{n,}s} :=\Vhn{s}{1}$,
$\U_{h{,n}} := \Vhn{s}{1} \cap H_0^1(\Omega_s;\Gamma_s^D)$,
$\P_{h{,n}} := \Vhn{f}{1}$ and
$\X_{h{,n}} := \V_{h,{n,}f} \times \P_{h{,n}} \times \V_{h,{n,}s} \times \U_{h{,n}}$. For the fluid part, the pair ($\V_{h,{n,}f}, \P_{h{,n}}$) is the inf-sup stable Taylor-Hood finite element space of lowest order; see e.g.~\cite{GiRa1986}.

The interface conditions are then imposed weakly using Nitsche's method,
see \eqref{pb2_4}--\eqref{pb2_5}, and so-called ghost penalty terms \eqref{pb2_6}, \eqref{pb2_7}, \eqref{pb2_9} are added
around the interface zone. These
extend the coercivity of the bilinear form over the interface cells and
increase stability. The fully discrete weak formulation then reads as follows.
\begin{problem}
\label{problem_discrete}
For $n=1,\ldots,N$ find
the fluid velocity, pressure, solid velocity and displacement
$U_h^n:= (v_f^h, p^h, v_s^h, u^h) := (v_f^{h, n}, p^{h, n}, v_s^{h, n}, u^{h, n})$
in $(v_{f{,n}}^D, 0, 0, u_{n}^D) + \X_{h{,n}}$, where
$(v_f^{h, n-1}, p^{h, n-1}, v_s^{h, n-1}, u^{h, n-1})$
are the solutions of the previous time step, such that
for all $\Psi^h:=(\phi_f^h, \psi^h, \phi_s^h, \xi^h) \in \X_{h{,n}}$:
\begin{align*}
A(U_h^n, \Psi^h) := F(\Psi^h),
\end{align*}
where
\begin{align}
A(&U_h^n, \Psi^h) := \rho_f \prodf{v_f^h}{\phi_f^h} +
    \rho_f k\prodf{ v_f^h \cdot \nabla v_f^h}{\phi_f^h} +
    k \prodf{\sigma_f(v_f^h, p^h)}{\nabla \phi_f^h}   \label{pb2_1} \\[-.5pt]
  &+
    k \prodf{\nabla \cdot v_f^h}{\xi^h}- k\prodd{ \rho_f \nu_f (\nabla v_f^h)^\top n_f}{\phi_f^h}
    + \rho_s \prods{v_s^h}{\phi_s^h}
     \label{pb2_2} \\[-.5pt]
  &+ k \prods{\sigma_s^h}{\nabla \phi_s^h}+ \rho_s k\prods{ v_s^h \cdot \nabla v_s^h}{\phi_s^h}
    + \prods{u^h +k(v_s^h \cdot \nabla u^h - v_s^h)}{\psi^h}
    \label{pb2_3} \\[-.5pt]
  &+ \frac{k}{h} \rho_f \nu_f \gamma_N
    \prodi{v_f^h - v_s^h}{\phi_f^h - \phi_s^h}
    - k \prodi{\sigma_f{(v_f^h, p^h)} \cdot n_f}{\phi_f^h - \phi_s^h}
    \label{pb2_4} \\[-.5pt]
  &- k \prodi{v_f^h - v_s^h}{\sigma_f (\phi_f^h, -\xi^h) \cdot n_f}
    \label{pb2_5} \\[-.5pt]
  &+ 2 \rho_f \nu_f k g_{v_f}^{h,n,w}(v_f^h, \phi_f^h) +
    \rho_s g_{v_s}^{h,n,w}(v_s^h, \phi_s^h) \label{pb2_6} \\
    &+k g_p^{h,n,w}(p^h, \xi^h) +
    2 \mu_s k g_u^{h,n,w}(u^h, \phi_s^h), \label{pb2_7} \\[-.5pt]
    F(&\Psi^h)
  := \rho_f k \prodf{f}{\phi_f^h} +
    \rho_f \prodf{v_f^{h, n-1}}{\phi_f^h} +
    \rho_s k \prods{ f}{\phi_s^h}  \label{pb2_8} \\[-.5pt]
  &+
    \rho_s \prods{v_s^{h, n-1}}{\phi_s^h}+ \prods{u^{h, n-1}}{\psi^h}
    + \rho_s g_{v_s}^{h,n,w}(v_s^{h, n-1}, \phi_s^h), \label{pb2_9}
\end{align}
and $\sigma_s^h := \sigma_s(u^h)$.
As previously already mentioned, here $k=t_n - t_{n-1}>0$ is the (uniform) time step size, $h>0$ the spatial discretization
parameter, namely the maximum element size,
and $\gamma_N>0$ denotes the Nitsche parameter. Note that the notation $\sigma_f (\phi_f^h, -\xi^h)$
is defined in~\eqref{fstress} by applying the Cauchy stress operator to the test functions $\phi_f^h$ and $-\xi^h$.
The ghost penalty term in \eqref{pb2_9} results from the time discretization of the term $g_{v_s}^{h,n,w}(\partial_t v_s^h, \phi_s^h)$ which is required in the stability analysis of the system, which is ongoing work.
\end{problem}

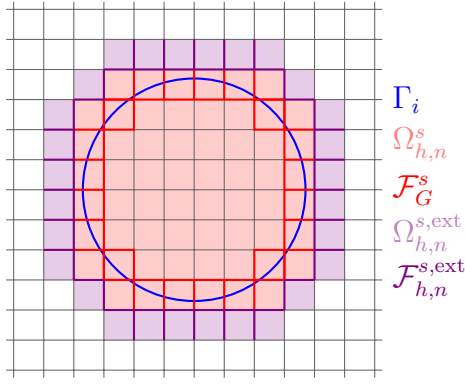
\begin{figure}[t]
  \centering
  \begin{tikzpicture}[font = \scriptsize, scale = 2.0]
    \filldraw[color=violet!20,fill=violet!20,thin] (-0.6,-1.0) rectangle (0.6,1.0);
    \filldraw[color=violet!20,fill=violet!20,thin] (-1.0,-0.6) rectangle (1.0,0.6);
    \filldraw[color=violet!20,fill=violet!20,thin] (-0.8,-0.8) rectangle (0.8,0.8);
    \filldraw[color=red!20,fill=red!20,thin] (-0.6,-0.8) rectangle (0.6,0.8);
    \filldraw[color=red!20,fill=red!20,thin] (-0.8,-0.6) rectangle (0.8,0.6);
    \draw[step=0.2,black!60,thin] (-1.25,-1.25) grid (1.25,1.25);
    \draw[color=blue,thick] (0.0,0.0) circle (0.74);

    \foreach \x in {0.0,0.2,...,0.8}
    {
      \draw[color=red,thick] (\x-0.4,0.6) -- (\x-0.4,0.8);
      \draw[color=red,thick] (\x-0.4,-0.6) -- (\x-0.4,-0.8);
      \draw[color=red,thick] (0.6,\x-0.4) -- (0.8,\x-0.4);
      \draw[color=red,thick] (-0.6,\x-0.4) -- (-0.8,\x-0.4);
      \draw[color=violet,thick] (\x-0.4,0.8) -- (\x-0.4,1.0);
      \draw[color=violet,thick] (\x-0.4,-0.8) -- (\x-0.4,-1.0);
      \draw[color=violet,thick] (0.8,\x-0.4) -- (1.0,\x-0.4);
      \draw[color=violet,thick] (-0.8,\x-0.4) -- (-1.0,\x-0.4);
    }

    \draw[color=red,thick] (-0.6,0.6) -- (0.6,0.6);
    \draw[color=red,thick] (0.6,-0.6) -- (0.6,0.6);
    \draw[color=red,thick] (-0.6,-0.6) -- (-0.6,0.6);
    \draw[color=red,thick] (-0.6,-0.6) -- (0.6,-0.6);
    \draw[color=red,thick] (-0.6,-0.4) rectangle (-0.4,-0.6);
    \draw[color=red,thick] (0.6,0.4) rectangle (0.4,0.6);
    \draw[color=red,thick] (0.6,-0.4) rectangle (0.4,-0.6);
    \draw[color=red,thick] (-0.6,0.4) rectangle (-0.4,0.6);

    \draw[color=violet,thick] (-0.6,0.8) -- (0.6,0.8);
    \draw[color=violet,thick] (-0.6,-0.8) -- (0.6,-0.8);
    \draw[color=violet,thick] (0.8,-0.6) -- (0.8,0.6);
    \draw[color=violet,thick] (-0.8,-0.6) -- (-0.8,0.6);

    \draw[color=violet,thick] (-0.8,0.6) -- (-0.6,0.6);
    \draw[color=violet,thick] (0.8,0.6) -- (0.6,0.6);
    \draw[color=violet,thick] (-0.8,-0.6) -- (-0.6,-0.6);
    \draw[color=violet,thick] (0.8,-0.6) -- (0.6,-0.6);

    \draw[color=violet,thick] (0.6,0.6) -- (0.6,0.8);
    \draw[color=violet,thick] (-0.6,0.6) -- (-0.6,0.8);
    \draw[color=violet,thick] (0.6,-0.6) -- (0.6,-0.8);
    \draw[color=violet,thick] (-0.6,-0.6) -- (-0.6,-0.8);

    \draw (1.25,0.6) node[anchor=west] {\large \textcolor{blue}{$\Gamma_i$}};
    \draw (1.25,0.3) node[anchor=west] {\large \textcolor{red!50}{$\Ocomp{s}$}};
    \draw (1.25,0.0) node[anchor=west] {\large \textcolor{red}{$\mathcal{F}_G^s$}};
    \draw (1.25,-0.3) node[anchor=west] {\large \textcolor{violet!50}{$\Oext{s}$}};
    \draw (1.25,-0.6) node[anchor=west] {\large \textcolor{violet}{$\FExt{s}$}};
  \end{tikzpicture}
  \caption{Visualization of the sets $\mathcal{F}_G^s$ and $\FExt{s}$ used in the definitions of the ghost penalty functions \eqref{gp} and \eqref{gp_ext} for a circular solid domain.}
  \label{fig:gp_sets}
\end{figure}

\begin{remark}
Let us comment on the Nitsche coupling terms in \eqref{pb2_5}--\eqref{pb2_6}. The first term in~\eqref{pb2_5} ensures continuity of the velocities across the interface in a weak sense, the second term in~\eqref{pb2_5} leads to a consistent imposition of this condition. The term in~\eqref{pb2_6} has first been proposed in~\cite{BURMAN2014FSI}. It is symmetric with respect to the fluid functions $v_f^h$ and $\phi_f^h$, but skew-symmetric in the pressure variable ($p^h$ vs. $\xi^h$). This is important in the stability analysis in~\cite{BURMAN2014FSI,BurmanFernandezFrei2020}, as the pressure terms cancel each other out by ``diagonal testing'' (i.e., when choosing $\xi^h = p^h$).
\end{remark}
Let $G_h^{{n}}:=\defset{T \in \T_h}{T \cap \Gamma_i{(t_n)} \neq \0}$ be the set of interface cells. Moreover let $\mathcal{F}$ be the set of
all element faces of the triangulation $\T_h$ and
let $\FG{f}:=\FG{f}(t_n) \subseteq \mathcal{F}$ denote the set of faces
$F = \cl{T}_1 \cap \cl{T}_2$ of $\tria{l}$
that do not lie on the boundary $\bd\Omega$,
such that at least one of the cells $T_j$ is intersected by the interface,
i.e., $T_j \cap \Gamma_i \ne \0$ for $j \in \set{1, 2}$.
Analogously, we define $\FG{s}:=\FG{s}(t_n)$
as the set of corresponding faces of the triangulation $\tria{s}$, which is displayed in Figure \ref{fig:gp_sets}.
For a cell $T$ cut by the interface we denote by $T\tsb{in}$
the part of the cell inside the considered subdomain.
For arbitrary neighboring cells ${T}$, ${T}'$ with common face $F = {T} \cap {T}'$
and unit normal $n$
(not necessarily $n_f$ or $n_s$)
let $\jump{\partial^i_n \varphi} := \partial^i_n \varphi|_{{T}} - \partial^i_n \varphi|_{{T}'}$ denote the jump in the $i$-th normal
of a function $\varphi$, where the orientation of the normal is arbitrary.
Using
$\Jump{i} := \iprod[F]{\jumpn1i}{\jumpn2i}$,
the ghost penalty functions (with positive parameters
$\gamma_{v_f}$, $\gamma_p$, $\gamma_{v_s}$, $\gamma_u$)
are defined as follows:
\newcommand\sumFf{\sum_{F \in \FG{f}}}
\newcommand\sumFKf{\sumFf \sum_{{T}\mathpunct: F \in \cl{{T}}}}
\newcommand\sumFs{\sum_{F \in \FG{s}}}
\newcommand\sumFKs{\sumFs \sum_{{T}\mathpunct: F \in \cl{{T}}}}
\begin{align} \label{gp}
  \begin{split}
  \ghw[v_f]{\varphi_1}{\varphi_2}
  &:= \gamma_{v_f} \sumFKf
    w(\kappa_{{T}}) \Bigl( h \Jump1 + \frac{h^3}{4} \Jump2 \Bigr), \\
  \ghw[p]{\varphi_1}{\varphi_2}
  &:= \gamma_p \sumFKf
  w(\kappa_{{T}}) h^3 \Jump1, \\
  \ghw[v_s]{\varphi_1}{\varphi_2}
  &:= \gamma_{v_s} \sumFKs w(\kappa_{{T}}) h^3 \Jump1, \\
  \ghw[u]{\varphi_1}{\varphi_2}
  &:= \gamma_u \sumFKs w(\kappa_{{T}}) h \Jump1.
  \end{split}
\end{align}
Here we apply a novel weight function \cite{FrKnStWeWi24_ENUMATH}
\begin{align*}
  w\colon [0,1] \to \Bigl[ \frac12 \wmax^{-1}, \frac12 \wmax \Bigr],
  \quad \kappa_{{T}} \mapsto \frac12 \wmax^{1 - 2 \kappa_T},
\end{align*}
with $\wmax \ge 1$,
which scales the ghost penalties depending on the cell cuts by taking
the portion of the inside cell part,
$\kappa_T := \text{meas}(T_\text{in}) / \text{meas}(T)$, as the argument.
Since $w(\frac12) = \frac12$ for any $w_{\max} \ge 1$,
an equally cut cell delivers the same ghost penalty contribution
as in the conventional terms, but if the inside part of the cell gets smaller,
so that the approximation gets worse, this contribution becomes larger and approaches $w_{\max}$ if
${T}_\text{in} \to \0$.
Thus we penalize ``bad cuts'' more severely
while ``good cuts'' (where a sufficiently large portion of the cell lies inside)
are penalized less severely.
Moreover,
  the conventional ghost penalty terms are recovered as the special case
  where $\wmax = 1$, hence $w \equiv \frac12$.

\begin{figure}
  \centering
  \def\plot#1{\includegraphics[width=0.48\textwidth,
  trim=0cm 25cm 20cm 0cm,clip]{pics/enumath_v_f_w_#1.png}}
  \plot1
  \plot3
  \includegraphics[width=\textwidth,
  trim=0cm 19cm 0cm 12cm,clip]{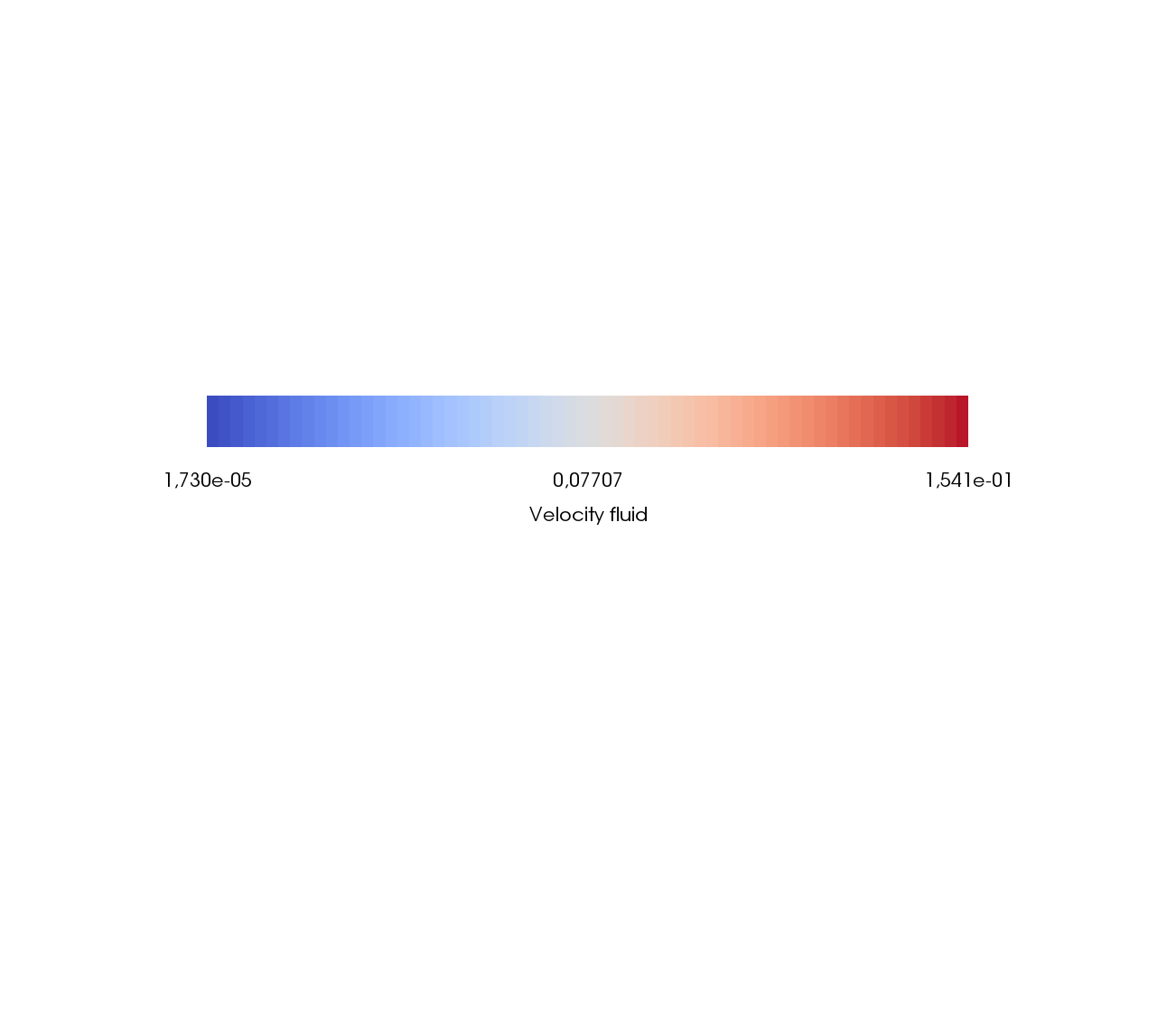}
  \caption{Euclidean norm of the fluid velocity $\|v_f\|_2$ on cut cells using $\wmax=1$ (left) and $\wmax=3$ (right). The values are taken from the ''modified flow around a cylinder`` computation in \cite{FrKnStWeWi24_ENUMATH}.}
  \label{fig:weights}
\end{figure}
A larger value of $\wmax$ reduces the fluctuation of the solution on cut cells, as can be observed in~\cref{fig:weights}.

\subsection{Stabilization of the convection terms}
For large interface movements
the structure convection terms have to be stabilized.
Therefore, we apply the Streamline Upwind Petrov–Galerkin (SUPG)~\cite{BrHu82} method
by adding the artificial diffusion terms
\begin{align*}
 S_{v_s}(U_h^n)(\Psi^h)&:=\delta_{h,k}^{v_s} \prods{v_s^h \cdot \nabla v_s^h}{v_s^h \cdot \nabla \phi_s^h} \quad \text{and} \\
 S_u(U_h^n)(\Psi^h)&:=\delta_{h,k}^u \prods{v_s^h \cdot \nabla u^h}{v_s^h \cdot \nabla \psi^h},
\end{align*}
where
\begin{align*}
  \delta_{h,k}^{v_s / u}
  :=
  \delta_0^{v_s / u} \frac{h^2}
  {6 \mu_s/\rho_s + h \norm{v_s^{h,j}}_{\Omega_s} + h/k}
  \quad \text{with} \quad \delta_0^{v_s} > 0, \ \delta_0^u > 0,
\end{align*}
that can be extended to a consistent formulation as follows:
\begin{align*}
 S_\text{supg}^{v_s}(U_h^n)(\Psi^h)&:=\delta_{h,k}^{v_s} \prods{\rho_s (v_s^h-v_s^{h,n-1})+k\rho_s v_s^h \cdot \nabla v_s^h-k\nabla \cdot \sigma_s^h-k\rho_s f}{v_s^h \cdot \nabla \phi_s^h}, \\
 S_\text{supg}^u(U_h^n)(\Psi^h)&:=\delta_{h,k}^u \prods{u^h-u^{h,n-1}+k v_s^h \cdot \nabla u^h-k v_s^h}{v_s^h \cdot \nabla \psi^h}, \text{ and} \\
 S_\text{supg}(U_h^n)(\Psi^h) &:= S_\text{supg}^{v_s}(U_h^n)(\Psi^h) + S_\text{supg}^{u}(U_h^n)(\Psi^h).
\end{align*}

%%%%%%%%%%%%%%%%%%%%%%%%%%%%%%%%%%%%%%%%%%%%%%%%%%%%%%%%%%%%%%%%%%
\subsection{Implicit Extensions to deal with the Moving Interface}
\label{subsec.ext}

In~\eqref{pb2_8}--\eqref{pb2_9} the quantities $v_s^{h,n-1}$ and $u^{h,n-1}$ are needed on $\Omega_s(t_n)$, but are defined on $\Omega_{s,h}^{n-1}$ (similarly for $v_f^{h,n-1}$). To handle this mismatch, we extend $v_s^{h,n-1}$ already in the previous time step ($t_{n-1}$) smoothly by means of ghost penalty terms, following~\cite{LehOl19}. For simplicity, we assume that the time step $k$ is chosen small enough such that
\begin{align*}
  k\leq v_{i,max} h,
\end{align*}
where $v_{i,max}$ is the maximum velocity of the interface displacement. This means that the interface moves by at most one layer of mesh cells per time step. Thus, it suffices to extend the computational mesh at time $t_n$ by one layer:
\begin{align*}
 \triaext{l} &:= \{T \in \T_h \ | \ \exists T' \in \tria{l} : \cl T \cap \cl{T'} \neq \emptyset\} \supseteq \tria{l} \quad \text{and} \\
 \Oext{l} &:= \bigcup_{T \in \triaext{l}} \cl T \supseteq \Ocomp{l} \quad \text{for } l \in \{f,s\}.
\end{align*}
We compute the extensions of the solution, denoted by $U_h^{n,\text{ext}}=(v_{f,\text{ext}}^h,p\tsb{ext}^h,v_{s,\text{ext}}^h,u\tsb{ext}^h)$,
by applying additional ghost penalty terms that ensure a smooth continuation on the extension cells:
\newcommand\sumFfext{\sum_{F \in \FExt{f}}}
\newcommand\sumFKfext{\sumFfext \sum_{{T}\mathpunct: F \in \cl{{T}}}}
\newcommand\sumFsext{\sum_{F \in \FExt{s}}}
\newcommand\sumFKsext{\sumFsext \sum_{{T}\mathpunct: F \in \cl{{T}}}}
\begin{align} \label{gp_ext}
  \begin{split}
  \ghwext[v_f]{\varphi_1}{\varphi_2}
  &:= \gammaext{v_f} \sumFKfext
    w(\kappa_{{T}}) \Bigl( h \Jump1 + \frac{h^3}{4} \Jump2 \Bigr), \\
  g_{p,\text{ext}}^{h,n,w}(\varphi_1, \varphi_2)
  &:= \gammaext{p} \sumFKfext
  w(\kappa_{{T}}) h^3 \Jump1, \\
  \ghwext[v_s]{\varphi_1}{\varphi_2}
  &:= \gammaext{v_s} \sumFKsext w(\kappa_{{T}}) h^3 \Jump1, \\
  \ghwext[u]{\varphi_1}{\varphi_2}
  &:= \gammaext{u} \sumFKsext w(\kappa_{{T}}) h \Jump1, \\
  \ghwext[u,\psi]{\varphi_1}{\varphi_2}
  &:= \gammaext{u,\psi} \sumFKsext w(\kappa_{{T}}) h \Jump1%end akw
  \end{split}
\end{align}
where $\FExt{l}:=\defset{F \in \mathcal{F} \setminus \bd\Oext{l}}{\exists T \in \triaext{l} \setminus \tria{l} : F \in \cl T} \setminus \FG{l}(t_n)$
for $l \in \{f,s\}$ (see Figure \ref{fig:gp_sets}) and with positive parameters \gammaext{v_f}, \gammaext{p}, \gammaext{v_s}, \gammaext{u} and \gammaext{u,\psi}. Compared to the ghost penalties, the additional term \gammaext{u,\psi} serves to extend the test function $\psi$ onto the additional cell layer.

The resulting modified fully discrete weak form is then given as follows.
\begin{problem}
\label{problem_discrete_ext}
Find $U_h^{n,\textup{ext}}=(v_{f,\textup{ext}}^h,p\tsb{ext}^h,v_{s,\textup{ext}}^h,u\tsb{ext}^h)$
for the time steps $n = 1, \dots, N$, such that
\begin{align*}
 A\tsb{ext}(U_h^{n,\textup{ext}})(\Psi^h):={}&A(U_h^{n,\textup{ext}})(\Psi^h)+ S_\textup{supg}(U_h^n)(\Psi^h) +2 \rho_f \nu_f k \ghwext[v_f]{v_{f,\textup{ext}}^h}{\phi_f^h}\\
 &+k \ghwext[p]{p\tsb{ext}^h}{\xi^h}
 +\rho_s k \ghwext[v_s]{v_{s,\textup{ext}}^h}{\phi_s^h} \\
 &+2 \mu_s k \ghwext[u]{u\tsb{ext}^h}{\phi_s^h} + k \ghwext[u,\psi]{u\tsb{ext}^h}{\psi^h}\\
 &=F(\Psi^h),
\end{align*}
where $A(\fcdot)(\fcdot)$ and $F(\fcdot)$ are defined in Problem~\ref{problem_discrete}.
\end{problem}

To define a level set function $\Phi$, which describes the position of the interface, we use the displacement $u\tsb{ext}^h(x,t_n)$ in the extended solid domain $\Oext{s}$ to set
\begin{align}
  \Phi(x,t_n)\label{Phin}
  &:= \begin{cases}
    \Phi^{{0}}(x - u\tsb{ext}^h(x,t_n)), & x \in \Oext{s}, \\
    1, & x \notin \Oext{s},
  \end{cases}
  %\\
\end{align}
where $\Phi^0\colon \Omega \to \R$ describes the initial configuration of the solid subdomain $\Omega_s(0)$ by satisfying the following equivalences:
\begin{align*}
 {\Phi^0(x)}&<0 \iff x \in \Omega_{s}(0), \\
 {\Phi^0(x)}&=0 \iff x \in \Gamma_i(0), \\
 {\Phi^0(x)}&>0 \iff x \notin \Omega_{s}(0) \cup \Gamma_i(0).
\end{align*}
As we assume that the interface $\Gamma_i$ does not leave $\Oext{s}$ within the $n$-th time step, we can safely set $\Phi=1$ in points $x \notin \Oext{s}$ in~\eqref{Phin}.
\begin{remark}
For better readability we waive the additional ``ext''-index and denote the extended solution simply by $U_h^{n}=(v_{f}^h,p^h,v_{s}^h,u^h)$ in the following.
\end{remark}

\section{Fluid-Structure-Contact Interaction}
\label{sec_cont}
We extend the model towards the situation that the solid $\Omega_s(t)$ can come into contact with a planar wall $\Gamma_w\subset\partial\Omega$. For simplicity, we assume that the wall $\Gamma_w$ is placed at the lower boundary of $\Omega$.

\subsection{Relaxation of the Contact Condition}

Modelling the transition from a fluid-structure interaction problem to full solid-solid contact introduces a number of severe numerical challenges, such as a non-smooth switch between fluid-structure interaction and solid-solid contact conditions, a topology change in the fluid domain $\Omega_f(t)$ when the fluid between the solids vanishes, and the Navier-Stokes no-collision paradox. The latter states that contact is not possible when the incompressible Navier-Stokes equations are imposed in the fluid in combination with no-slip conditions on the boundary $\Gamma_w$ and/or no-slip interface conditions on $\Gamma_i(t)$~\cite{Hillairet2d, HeslaPhD, HillairetTakahashi3d}.

It is also questionable whether such a model would accurately describe the underlying physics, as due to surface roughness of the solids some seepage of the fluid through the contact surface should always be possible~\cite{AgerWalletal, BurmanFernandezFreiGerosa2024}. Therefore, a number of works on FSI-contact have recently considered to model a porous fluid layer between the contacting surfaces~\cite{BurmanFernandezFrei2020, AgerWalletal, BurmanFernandezFreiGerosa2024, GerosaMarsden2024}. The simplest possibility is to use the Navier-Stokes equations in the small fluid layer between the contacting solids as well. This can be achieved by a relaxation of the contact conditions.

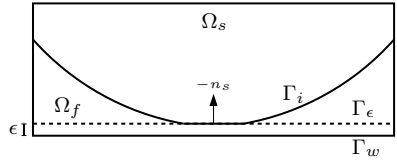
\begin{figure}[t]
\centering
\begin{picture}(0,0)%
\includegraphics{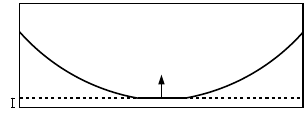}%
\end{picture}%
\setlength{\unitlength}{1657sp}%
\begingroup\makeatletter\ifx\SetFigFont\undefined%
\gdef\SetFigFont#1#2{%
  \fontsize{#1}{#2pt}%
  \selectfont}%
\fi\endgroup%
\begin{picture}(5808,2383)(1426,-3311)
\put(4321,-1321){\makebox(0,0)[lb]{\smash{{\SetFigFont{8}{9.6}{\color[rgb]{0,0,0}$\Omega_s$}%
}}}}
\put(2116,-2581){\makebox(0,0)[lb]{\smash{{\SetFigFont{8}{9.6}{\color[rgb]{0,0,0}$\Omega_f$}%
}}}}
\put(1441,-2919){\makebox(0,0)[lb]{\smash{{\SetFigFont{8}{9.6}{\color[rgb]{0,0,0}$\epsilon$}%
}}}}
\put(6571,-2626){\makebox(0,0)[lb]{\smash{{\SetFigFont{8}{9.6}{\color[rgb]{0,0,0}$\Gamma_{\epsilon}$}%
}}}}
\put(6571,-3211){\makebox(0,0)[lb]{\smash{{\SetFigFont{8}{9.6}{\color[rgb]{0,0,0}$\Gamma_w$}%
}}}}
\put(4231,-2221){\makebox(0,0)[lb]{\smash{{\SetFigFont{5}{6.0}{\color[rgb]{0,0,0}$-n_s$}%
}}}}
\put(5536,-2401){\makebox(0,0)[lb]{\smash{{\SetFigFont{8}{9.6}{\color[rgb]{0,0,0}$\Gamma_i$}%
}}}}
\end{picture}%
\caption{Relaxation of the contact condition with a planar wall: We impose the no-penetration condition already at a distance $\epsilon>0$ from the lower wall.
\label{fig:relaxedcont}
}
\end{figure}

Therefore, let $g_0^n(x)$ be the current distance of a point $x\in \Gamma_i(t_n)$ to the lower (planar) wall $\Gamma_w$. The classical no-penetration condition (see e.g., \cite{Wriggers2006}) reads in an Eulerian formulation
\begin{align}
(u^{h,n}-u^{h,n-1}) \cdot n_{{s}} \le g_0^n, \label{cont1}
\end{align}
where $n_{{s}}$ is the normal of the solid subdomain on the interface $\Gamma_i$, which is equal to the normal of the wall $\Gamma_w$ in the case of contact.
To keep a fluid layer for all times, we relax~\eqref{cont1} by a small $\epsilon>0$ and impose that a distance of $\epsilon$ will always remain between solid and wall:
\begin{align}
(u^{h,n}-u^{h,n-1}) \cdot n_{{s}} \le g_\epsilon^n := g_0^n-\epsilon.\label{conteps}
\end{align}
The configuration is illustrated in Figure~\ref{fig:relaxedcont}.

 \subsection{Combination of FSI and contact conditions}

When it comes to contact, the dynamic FSI condition
$\sigma_f {\cdot} n_s = \sigma_s {\cdot} n_s$
is not valid any more, due to an additional contact force.
Following~\cite{BurmanFernandezFrei2020},
we model this by a means of a Lagrange multiplier $\lambda\leq 0$,
which acts in direction $-n_{{s}}$ (see Figure~\ref{fig:relaxedcont})
when augmenting the discrete formulation,
i.e., Problem~\ref{problem_discrete_ext}:
\begin{align}\label{contLagr}
A\tsb{ext}(U_h^n)(\Psi^h) -k\prodi{\lambda n_s}{\phi_s^h} = F(\Psi^h).
\end{align}
By collecting all terms in~\eqref{contLagr} that involve the test function $\phi_s^h$ on the interface $\Gamma_i$, we see that this corresponds to the interface condition
\begin{align} \label{concon_lam}
 &n_{{s}}^\top \cdot \left(\frac{\rho_f \nu_f \gamma_N}{h} (v_f^h-v_s^h) - \sigma_s^h \cdot n_{{s}}+\sigma_f^h \cdot n_{{s}} \right) - \lambda=0 \quad \text{ on } \Gamma_i,
\end{align}
where $\sigma_f^h := \sigma_f(v_f^h, p^h)$ and $\sigma_s^h:= \sigma_s(u^h)$.
In the case $\lambda=0$, we recover the FSI interface condition, for
$\lambda<0$ an additional contact force is added. As the latter is only
allowed when contact is active (i.e., $(u^h-u^{h,n-1}) \cdot n_s = g_\epsilon^n$),
the following conditions have to be satisfied:
\begin{align*}
 &(u^h-u^{h,n-1}) \cdot n_{{s}} \le g_\epsilon^n, \quad \lambda \le 0, \quad \Big((u^h-u^{h,n-1}) \cdot n_{{s}}-g_\epsilon^n \Big)\lambda=0.
 \end{align*}
 Following Alart \& Curnier~\cite{AlartCurnier1991}, these three conditions are equivalent to the following (non-smooth) equality condition for arbitrary $\gamma_C>0$:
 \begin{align*}
 &\lambda=-\gamma_C \left[(u^h-u^{h,n-1}) \cdot n_{{s}}-g_0^n-\frac{\lambda}{\gamma_C} \right]_+ =: {-}\gamma_C \left[P_{\gamma_C}(\lambda,u^h) \right]_+,
\end{align*}
where the operator $[f]_+:=\max\{f,0\}$ stands for the positive part of a function $f$.
In the numerical implementation, we set
$\gamma_C:=\gamma_C^0 \mu_s h^{-1}$, $\gamma_C^0>0$, see~\cite{BurmanFernandezFrei2020}.

This condition can be inserted in the weak form by replacing $\lambda$ in (\ref{contLagr})
by $-\gamma_C \left[P_{\gamma_C}(\lambda,u^h) \right]_+$ and, according to (\ref{concon_lam}), by replacing $\lambda$ in $P_{\gamma_C}(\lambda,u^h)$ by the jump of discrete fluxes
\begin{align*}
 [\tilde \sigma_n]:=n_{{s}}^\top \cdot \left((\sigma_s^h-\sigma_f^h) \cdot n_f+\frac{\rho_f \nu_f \gamma_N}{h} (v_f^h-v_s^h) \right).
\end{align*}
Finally, the weak formulation with contact condition reads:
\begin{problem}
Find $U_h^n$ for the time steps $n = 1, \dots, N$, such that
\begin{align}\label{contPgamma}
 A\tsb{ext}(U_h^n)(\Psi^h)+\gamma_C k \prodi{\left[P_{\gamma_C}(U_h^n) \right]_+}{\phi_s^h n_{{s}}}=F(\Psi^h),
\end{align}
where
\begin{multline} \label{PgammaC}
 P_{\gamma_C}(U_h^n):=P_{\gamma_C}([\tilde \sigma_n],u^h) \\
 =(u^h-u^{h,n-1}) \cdot n_{{s}}-g_\epsilon^n-\gamma_C^{-1} n_{{s}}^\top \cdot \left((\sigma_s^h-\sigma_f^h) \cdot n_{{f}}+\frac{\rho_f \nu_f \gamma_N}{h} (v_f^h-v_s^h) \right).
\end{multline}
\end{problem}
\begin{remark}
We note that the switch between fluid-structure-interaction and contact condition is now included implicitly in the variational formulation~\eqref{contPgamma}, depending on the sign of the term $P_{\gamma_C}(U_h^n)$.
\end{remark}

%%%%%%%%%%%%%%%%%%%%%%%%%%%%%%%%%%%%%%%%%%%%%%%%%%%%%%%
\section{Nonlinear Solution by a Semi-Smooth Newton Approach}
\label{sec_sol}

In this section, we provide details for our nonlinear solution algorithm.
The alorithm is based on Newton's method. Specifically, the Jacobian
is derived analytically by computing the Fr\'echet derivative by hand.

Let us start by deriving Newton's method for Problem~\ref{problem_discrete_ext}
without contact term.
With the step length $\alpha^j \in (0,1]$
determined by a backtracking line search,
Newton's method takes the following form.
\begin{problem}
\label{problem_Newton}
Let Problem~\ref{problem_discrete_ext} be given.
At any given time point $t_n$ for $n=1,\ldots, N$
solve the nonlinear problem by a Newton iteration indexed by $j$, as follows.
Given an initial guess $U_h^{n,0} \in (v_{f{,n}}^D, 0, 0, u_{n}^D) + \X_{h{,n}}$,
such as $U_h^{n,0} := U_n^{n-1}$, find $\delta U_h \in \X_{h{,n}}$
for $j = 0, 1, 2, \dots$, such that for all $\Psi^h \in \X_{h{,n}}$:
\begin{align*}
  A\tsb{ext}'(U_h^{n,j})(\delta U_h,\Psi^h)
  &= -A\tsb{ext}(U_h^{n,j})(\Psi^h) + F(\Psi^h), \\
  U_h^{n,j+1} &= U_h^{n,j} + \alpha^j \delta U_h.
\end{align*}
Set $U_h^n := U_h^{n,j^*}$, where $j^*$ denotes the index at which
Newton's method converged. Here, a residual-based stopping criterion is used.
\end{problem}
In Problem~\ref{problem_Newton}, we need the derivative
\begin{align*}
  A&\tsb{ext}'(U_h^{n,j})(\delta U_h, \Psi^h) =
      \lim_{\epsilon \to 0} \frac{1}{\epsilon} \left(
      A\tsb{ext}(U_h^{n,j} + \epsilon \delta U_h)(\Psi^h) - A\tsb{ext}(U_h^{n,j})(\Psi^h) \right) \\
    &= \rho_f \prodf{\delta v_f^h}{\phi_f^h}
%%% Convection term
      + \rho_f k \prodf{\delta v_f^h \cdot \nabla v_f^{h,j} +
      v_f^{h,j} \cdot \nabla \delta v_f^h}{\phi_f^h}\\
%%% Stress fluid
    &+ k \prodf{\sigma_f^h (\delta v_f^h, \delta p^h)}{\nabla \phi_f^h}
      + k \prodf{\nabla \cdot \delta v_f^h}{\xi^h}
      - k \rho_f \nu_f
      \smash{\prodd{\nabla {\delta v_f^h}^\top n_f}{\phi_f^h}} \\
    &+ \rho_s \prods{\delta v_s^h}{\phi_s^h}
%%% Convection term
      + \rho_s k \prods{\delta v_s^h \cdot \nabla v_s^{h,j} +
      v_s^{h,j} \cdot \nabla \delta v_s^h}{\phi_s^h} \\
    &+ k \prods{{\sigma_s^h}'(u^{h,j})(\delta u^h)}{\nabla \phi_s^h}
      + %\frac{\gamma_{uv}}{h^3}
      \prods{\delta u^h + k(\delta v_s^h \cdot \nabla u^{h,j} +
      v_s^{h,j} \cdot \nabla \delta u^h - \delta v_s^h)}{\psi^h} \\
    &+ k h^{-1} \rho_f \nu_f \gamma_N
      \prodi{\delta v_f^h - \delta v_s^h}{\phi_f^h - \phi_s^h}
      - k \prodi{\sigma_f^h (\delta v_f^h, \delta p^h) \cdot n_f}
      {\phi_f^h - \phi_s^h} \\
    &- k \prodi{\delta v_f^h - \delta v_s^h}
      {\sigma_f^h (\phi_f^h, -\xi^h) \cdot n_f} \\
    &+ 2 \rho_f \nu_f k g_{v_f}^h (\delta v_f^h, \phi_f^h)
      + \rho_s g_{v_s}^h (\delta v_s^h, \phi_s^h)
      + k g_p^h (\delta p^h, \xi^h)
      + 2 \mu_s k g_u^h (\delta u^h, \phi_s^h)\\
       &+ 2 \rho_f \nu_f k \ghwext[v_f]{\delta v_f^h}{\phi_f^h}
 +k \ghwext[p]{\delta p^h}{\xi^h}
 +\rho_s k \ghwext[v_s]{\delta v_s^h}{\phi_s^h}\\
 &+ 2 \mu_s k \ghwext[u]{\delta u^h}{\phi_s^h}
 + k\ghwext[u,\psi]{\delta u^h}{\psi^h}
 +S_\text{supg}'(U_h^n)(\Psi^h)
\end{align*}
where
\begin{align*}
  {\sigma_s^h}'(u^{h,j})(\delta u^h)
  &:= 2 \mu_s E_s'(u^{h,j})(\delta u^h) +
    \lambda_s \tr \bigl(E_s'(u^{h,j})(\delta u^h)\bigr), \\[-1pt]
  E_s'(u^{h,j})(\delta u^h)
  &:= \frac12 \bigl( \nabla \delta u^h + (\nabla \delta u^h)^\top
    + (\nabla \delta u^h)^\top \cdot \nabla u^{j,h}
    + (\nabla u^{j,h})^\top \cdot \nabla \delta u^h \bigr).
\end{align*}
The derivative of the SUPG convection stabilization is
$S_\text{supg}' = {S_\text{supg}^{v_s}}' + {S_\text{supg}^{u}}' $, where
\begin{align*}
 {S_\text{supg}^{v_s}}'&(U_h^{n,j})(\delta U_h,\Psi^h)=\delta_{h,k}^{v_s} \Bigl( \prods{\rho_s \delta v_s^h+k \rho_s (\delta v_s^h \cdot \nabla v_s^{h,j}+v_s^{h,j} \cdot \nabla \delta v_s^h)\\
 &-k\nabla \cdot {\sigma_s^h}'(u^{h,j})(\delta u^h)}{v_s^{h,j} \cdot \nabla \phi_s^h}
 -\rho_s \prods{v_s^{h,n-1}+kf}{v_s^{h,j} \cdot \nabla \phi_s^h} \\
 &+\prods{\delta v_s^h}{((\rho_s (v_s^{h,j}-v_s^{h,n-1})+k\rho_s v_s^{h,j} \cdot \nabla v_s^{h,j}-k\nabla \cdot \sigma_s^{h,j}-k\rho_s f) \cdot \nabla \phi_s^h} \Bigr), \\
  {S_\text{supg}^u}'&(U_h^{n,j})(\delta U_h,\Psi^h)={}\\
  &\delta_{h,k}^u \Bigl( \prods{\delta u^h-u^{h,n-1} + k(\delta v_s^h \cdot \nabla u^{h,j}+v_s^{h,j} \cdot \nabla \delta u^h - \delta v_s^h)}{v_s^{h,j} \cdot \nabla \psi^h} \\
 &\quad+\prods{\delta v_s^h}{(u^{h,j}-u^{h,n-1} + k(v_s^{h,j} \cdot \nabla u^{h,j} - v_s^{h,j})) \cdot \nabla \psi^h} \Bigr).
\end{align*}
We note that (shape) derivatives with respect to the sub-domains $\Omega_f(t_n)$ and $\Omega_s(t_n)$ (that also depend on the displacement $u^h$) are neglected. In this sense, our approach can be seen as a
simplified Newton method.

Finally, we need to take care of the derivation of the contact term
\begin{align*}
C(U_h^n, \Psi^h):= \prodi{\left[P_{\gamma_C}(U_h^n) \right]_+}{\phi_s^h n_{{s}}},
\end{align*}
with $P_{\gamma_C}(U_h^n)$ defined in~\eqref{PgammaC}. The maximum operator appearing in $[\fcdot]$ is not differentiable. Hence, we use a semi-smooth Newton method by setting
\begin{align*}
\tilde{C}'(U_h^n, \Psi^h) = \begin{cases} \prodi{P_{\gamma_C}'(U_h^n)(\delta U_h^n)}{\phi_s^h n_{{s}}}, &\text{ if } P_{\gamma_C}(U_h^n)>0, \\
0 &\text{else}.
\end{cases}
\end{align*}
Altogether, the semi-smooth Newton step reads as follows.

\begin{problem}
\label{problem_Newton_semi_smooth}
In extension to Problem~\ref{problem_Newton}:
\begin{algorithmic}
 \State Solve $A\tsb{ext}'(U_h^{n,j})(\delta U_h,\Psi^h) + S_\text{supg}'(U_h^{n,j})(\delta U_h, \Psi^h) + \tilde{C}'(U_h^{n,j})(\delta U_h,\Psi^h)$ for $\delta U_h$.
 \State $\qquad \quad \, =-A\tsb{ext}(U_h^{n,j})(\Psi^h)  -S_\text{supg}(U_h^{n,j})(\Psi^h) - C(U_h^j)(\Psi^h) + F(\Psi^h)$% + F_C(\Psi^h)$
 \State $U_h^{n,j+1} = U_h^{n,j} + \alpha^j \delta U_h$.
\end{algorithmic}
Set $U_h^n := U_h^{n,j^*}$, where $j^*$ denotes the index at which
Newton's method converged. Here, a residual-based stopping criterion is used.
\end{problem}

%%%%%%%%%%%%%%%%%%%%%%%%%%%%%%%%%%%%%%%%%%%%%%%%%%%%%%%
\section{Numerical Example: Bouncing elastic ball}
\label{sec_tests}
In this section, we conduct numerical simulations for a falling and bouncing
elastic ball. As previously mentioned, our code builds upon
the open-source finite element library deal.II~\cite{dealii2019design,dealII95},
in particular step-85,
and general FSI routines and basic nonlinear solver routines are based upon~\cite{Wi13_fsi_with_deal}.

\subsection{Geometry}\label{sec.geo}
\begin{figure}[tp]
\centering
\def\h{3.9}
\def\r{1.1}
\resizebox{0.5\textwidth}{!}{%
\begin{tikzpicture}[scale = 0.7]
  \draw (-4,0) -- (4,0) node [below] at (-4,0) {$-0.04$};
  \draw (4,0) -- (4,8) node [below] at (4,0) {$0.04$};
  \draw (4,8) -- node [above] {$\Gamma\tsb{top}$} (-4,8) node [left] at (-4,8) {$0.08$};
  \draw (-4,8) -- (-4,0) node [left] at (-4,0) {$0$};
  \draw (0,\h+\r) circle[radius=\r];
  \draw[dashed] node [below] {$0$} (0,0) -- (0,\h) -- (-4,\h) node [left] at (-4, \h) {$h_0$};
  \draw[thick] (0,\h+\r) -- node [above right] {$r$} (\r/2,\h+\r*0.13397);%[below right] {$r = 0.011$};
  \fill[black] (0,\h+\r) circle(1pt);
\end{tikzpicture}
}
\caption{Initial configuration for the bouncing elastic ball.
  All units are given in meters.
}
\label{fig:initial_conf}
\end{figure}
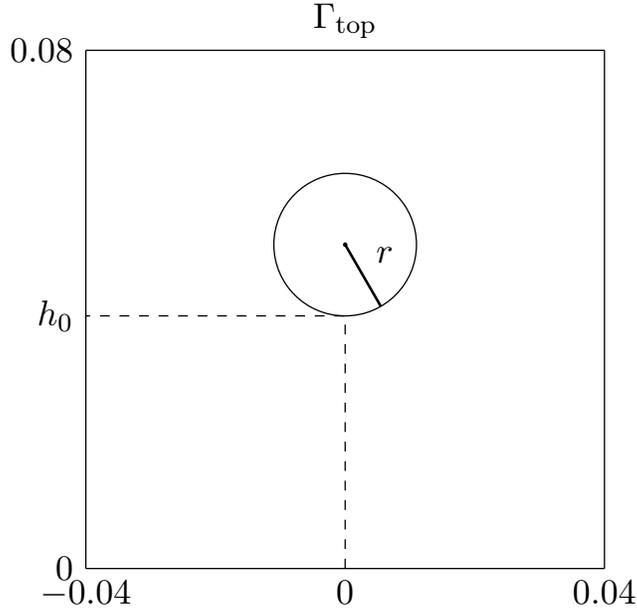

We consider a rectangular domain filled with an incompressible Newtonian fluid containing an elastic ball of radius $r = \SI{0.011}{m}$ and initial height $h_0 = \SI{0.039}{m}$, as seen in \cref{fig:initial_conf}.
The same computation is performed on two different types of meshes -- uniform quadrilateral meshes with element sizes of $\SI{1.25}{mm} \times \SI{1.25}{mm}$ with a total number of $37\, 660$ dofs after two global refinements, and a quadrilateral mesh of rectangular cells that decrease in height toward the bottom of the domain, giving
a minimum element size of $\SI{0.875}{mm} \times \SI{0.1875}{mm}$ on the bottom, a maximum element size of $\SI{1.875}{mm} \times \SI{1.625}{mm}$, and $55\, 786$ dofs.

\subsection{Boundary and initial conditions}
At the top boundary of the domain we employ a
slip boundary condition.
The vertical and bottom boundaries are no slip walls. Together, the boundary conditions read
\begin{align*}
 (\sigma_f n)_x = (v_f)_y &= 0 \text{ on } \Gamma\tsb{top},\\
 v_f &= 0 \text{ on } \bd{\Omega}\setminus{\Gamma\tsb{top}}.
\end{align*}
For the displacements $u$, homogeneous Neumann conditions are employed on the entire boundary $\bd{\Omega}$.
We employ zero initial conditions $v^0 = u^0 = 0$
for the velocities and displacements.
Both the fluid and the structure are subject to the
constant gravitational force field given as $f(x,y,t) = (0,-9.81)\,\si{m/s^2}$.

\subsection{Parameters}
The physical parameters of the fluid and the solid are given in \cref{tab:physprm}. For the ghost penalty parameters we choose $\gamma_{v_f} = 0.5$, $\gamma_{v_s} = \gamma_{p} = \gamma_u = 0.1$ with maximum weight $\wmax = 2.0$ and Nitsche parameter $\gamma_N = \num{1.0e+7}$. The extension parameters are prescribed by $\gammaext{u,\psi} = \num{1.0e-4}$ and $\gammaext{v_f} = \gammaext{p} = \gammaext{v_s} = \gammaext{u} = 10$, the SUPG-stabilization parameters by $\delta_0^{v_s} = \delta_0^u = \num{1.0e-5}$ and the contact parameter by $\gamma_C = 500$ with relaxed distance $\epsilon = \num{1.0e-4}$. An initial time step size of $k = \num{1.0e-4}$ is used at the start of the computation and the absolut tolerance for Newton's method is set to \num{1.0e-7}.

\begin{table}[tp]
  \caption{\label{tab:physprm}Physical parameters.}
  \centering
    \begin{tabular}{c@{\qquad}c@{\qquad}c@{\qquad}c@{\qquad}c}
    \toprule
    $\nu_f$ & $\rho_f$ & $\rho_s$ & $\mu_s$ & $\lambda_s$ \\
    \midrule
    \SI{7.0114e-5}{m^2/s} & \SI{1141}{kg/m^3} & \SI{1361}{kg/m^3} & \SI{20}{kPa} & \SI{80}{kPa}\\
    \bottomrule
  \end{tabular}
\end{table}

\subsection{Adaptive Time-stepping}
Extending the computational domain by one layer of cells for the interface
motion imposes a constraint on the magnitude of the deformation in a single
time step and consequently on the time step size. To reduce the computational
cost we adaptively refine and coarsen the time step size with respect to the
deformation. Additional refinement takes place in cases where the current time
step size does not yield convergence in Newton's method. Once the time
stepping is flagged for refinement, the time step size as well as the
  tolerance for Newton's method are reduced by a factor $\alpha_{{k}} =
0.1$. If no refinements have occurred during the prior five time steps, the
system is reset to five time steps earlier. Otherwise, it is reset to the
previous time. This approach makes it possible to continue the computation at
a point in time when errors due to too large time step sizes have not yet been
further amplified.
Both the time step size and Newton's tolerance
are later increased by a factor $\alpha_{k}^{-1}=10$,
at most up to the initial time step size,
after ten time steps have been taken without further refinement.
Coarsening after ten time steps (and not earlier)
serves to ensure that the solutions can be printed
on a uniform coarse time grid.

\subsection{Quantities of interest}
Motivated by the benchmark configuration~\cite{vonWahletal2021}
we define the following metrics to compare the simulations using the two different grids. Let $t^0>0$ be the time when the center of the ball reaches the height $h_0$ and $t^*$ be the time when the minimum distance between $\Gamma_i$ and $\Gamma_w$ is equal to the diameter of the ball, relative to the time $t^0$.
We denote by $v^*$ the average vertical velocity of the ball and by $f^*:= \int_{\Gamma_i} \sigma_f n e_2 \text{d}s$  the vertical force acting on the ball at time $t^0+t^*$ (where $e_2$ denotes the vertical unit vector). To quantify the fall and rebound of the ball, let $t\tsb{cont}$ be the time of the first contact between $\Gamma_i$ and $\Gamma_w$, and let $t\tsb{jump}$ be the time when the maximum
height $h\tsb{jump}$ of the ball is reached during the rebound, both relative to $t^0$. Furthermore, let $p\tsb{bc}$ be the pressure at the bottom center point $(0,0)$ and
\begin{align*}
 E\tsb{el}:=\int_{\Omega_s} \sigma_s : E_s \, \mathrm{d}x \quad \text{and} \quad E_{\text{kin},l}:=\frac{1}{2}\int_{\Omega_l} \rho_l \Vert v_l \Vert^2 \, \mathrm{d}x, \ l \in \{f,s\}
\end{align*}
be the elastic energy of the ball and the kinetic energies of the fluid and solid, respectively. Finally, we compare the average number of
Newton steps over all iterations.

\subsection{Discussion of our findings}
\subsubsection{Uniform mesh}
\begin{figure}
	\def\plot#1#2{%
		\begin{tikzpicture}
			\small
			\begin{axis}[
				xlabel={Time $t$},
				ylabel={Distance},
				xmin = #1, xmax=#2,
				ymin = 0,
				grid=both,
				width=185pt
				]
				\addplot[
				thin,
				color=blue
				]
				table[
				col sep=comma,
				header=false,
				x index=0,
				y index=1
				] {pics/dist_file_ref_0.csv};
				\addlegendentry{\small Level 0}
				\addplot[
				thin,
				color=red % Farbe der Punkte
				]
				table[
				col sep=comma,
				header=false,
				x index=0,
				y index=1
				] {pics/dist_file_ref_1.csv};
				\addlegendentry{\small Level 1}
				\addplot[
				thin,
%				color=red % Farbe der Punkte
				]
				table[
				col sep=comma,
				header=false,
				x index=0,
				y index=1
				] {pics/dist_file_ref_2.csv};
				\addlegendentry{\small Level 2}
				\addplot[
				thin,
				color=violet % Farbe der Punkte
				]
				table[
				col sep=comma,
				header=false,
				x index=0,
				y index=1
				] {pics/dist_file_ref_3.csv};
				\addlegendentry{\small Level 3}
			\end{axis}
		\end{tikzpicture}
	}
	\plot{0}{0.6}\hfill\plot{0.47}{0.6}
	\caption{Minimal distance between the interface and the bottom boundary against time using the uniform mesh.}
	\label{fig:plotsmesh1}
\end{figure}
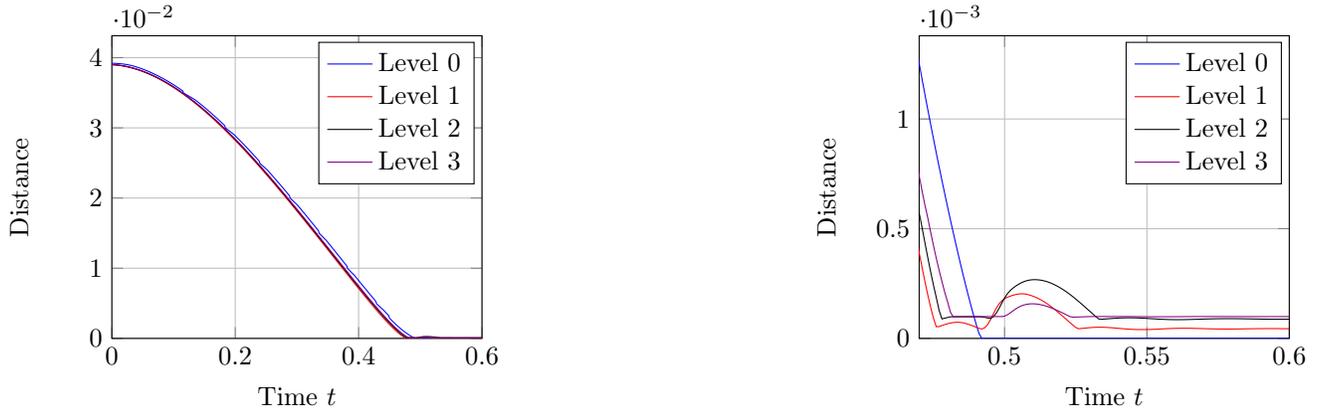

\begin{table}[tp]
	\caption{Quantities of interest on the uniform mesh on refinement level 0 to 3 (top to bottom). All quantities are given in SI units.}
	\centering
	\ifcase0
  \begin{tabular}{r*7c}
	\toprule
	DoFs & $t^0$ & $t^*$ & $v^*$ & $f^*$ & $t\tsb{cont}$ & $t\tsb{jump}$ & $h\tsb{jump}$ \\
	\midrule
	2695  & 0.2072 & 0.0623 &$-$0.1021 &$-$4.6861 & 0.3465 & ---    & ---    \\
	9928  & 0.2024 & 0.0613 &$-$0.1035 &$-$4.9170 & 0.2738 & 0.3037 & 0.000203 \\
	37660 & 0.2038 & 0.0613 &$-$0.1034 &$-$4.8917 & 0.2742 & 0.3068 & 0.000267 \\
	146648 & 0.2045 & 0.0615 &$-$0.1030 &$-$4.6639 & 0.2772 & 0.3053 & 0.000158 \\
	\bottomrule
\end{tabular}
$ $\\[1ex]
\begin{tabular}{r*6c}
	\bottomrule
	DoFs & $\max_t p\tsb{bc}$ & $\max_{Q_f} v_f$ & $\max_t E\tsb{el}$ & $\max_t E_{\text{kin},f}$ & $\max_t E_{\text{kin},s}$ & nNewton \\
	\midrule
	2695  & 677.3779  & 0.1485 & 0.0001 & 0.0087 & 0.0030 & 1.1952 \\
	9928  & 577.9210  & 0.1630 & 0.0031 & 0.0089 & 0.0032 & 1.2280 \\
	37660 & 575.0429  & 0.1705 & 0.0034 & 0.0090 & 0.0032 & 1.2496 \\
	146648 & 575.5481 & 0.1737 & 0.0021 & 0.0091 & 0.0031 & 1.3469\\
	\bottomrule
\end{tabular}
\or
\begin{tabular}{c@{\quad}c@{\quad}c@{\quad}c@{\quad}c@{\quad}c@{\quad}c@{\quad}c}
	\toprule
	DoFs & $t^0$ & $t^*$ & $v^*$ & $f^*$ & $t\tsb{cont}$ & $t\tsb{jump}$ & $h\tsb{jump}$\\
	\midrule
	2695  & 0.2072 & 0.0623 &$-$0.1021 &$-$4.6861 & 0.3465 & ---    & ---    \\
	9928  & 0.2024 & 0.0613 &$-$0.1035 &$-$4.9170 & 0.2738 & 0.3037 & 0.0002 \\
	146648 & 0.2045 & 0.0615 &$-$0.1030 &$-$4.6639 & 0.2772 & 0.3053 & 0.0002 \\
	\bottomrule\bottomrule
	DoFs & $\max_t p\tsb{bc}$ & $\max_{Q_f} v_f$ & $\max_t E\tsb{el}$ & $\max_t E_{\text{kin},f}$ & $\max_t E_{\text{kin},s}$ & \multicolumn{2}{c}{nNewton}\\
	\midrule
	2695  & 677.3779  & 0.1485 & 0.0001 & 0.0087 & 0.0030 & \multicolumn{2}{c}{1.1952} \\
	9928  & 577.9210  & 0.1630 & 0.0031 & 0.0089 & 0.0032 & \multicolumn{2}{c}{1.2280} \\
	37660 & 575.0429  & 0.1705 & 0.0034 & 0.0090 & 0.0032 & \multicolumn{2}{c}{1.2516} \\
	146648 & 575.5481 & 0.1737 & 0.0021 & 0.0091 & 0.0031 & \multicolumn{2}{c}{1.3469}\\
	\bottomrule
\end{tabular}
	\fi
	\label{tab:qoimesh}
\end{table}

First, we compare results on three different globally refined meshes; followed by comparing results between the uniform and the locally refined meshes. In addition to the mesh described in Section~\ref{sec.geo}, we use two coarser meshes (level 0 and 1) with a element size of $5$ resp.\,$2.5$mm and a finer mesh with a element size of $0.625$mm. The vertical distance between the ball and the bottom boundary for all four refinement levels is shown in \cref{fig:plotsmesh1}. We observe that no rebound occurs on the coarsest mesh. In fact, on refinement level 0 the contact condition is not penalized sufficiently and a part of the ball moves below the bottom boundary, indicating that the value of $\epsilon$ is too small for the coarse mesh size. 

All quantities of interest for the four refinement levels are provided in \cref{tab:qoimesh}. We observe significant differences between the coarsest mesh (level~0) and the other three mesh levels. The results on refinement levels 1--3 are in reasonable agreement. Comparing levels 1--3, we observe by compring $t^*$ that the impact happens later the finer the mesh is. This is due to the fact that the fluid forces slowing down the ball are better resolved on the fine mesh levels. Consequently, the impact velocity $v^*$ is smaller on the finest mesh (in terms of absolute values) and thus, the jump height $h_{jump}$ is reduced compared to level 2.

\subsubsection{Uniform and non-uniform mesh}

\begin{figure}
  \centering
  % trim={left bottom right top}
  \def\plot#1{\includegraphics[width=0.48\textwidth,
  trim=0cm 0cm 0cm 0cm,clip]{pics/solution_v_y_no_grid_#1.png}}
  \plot{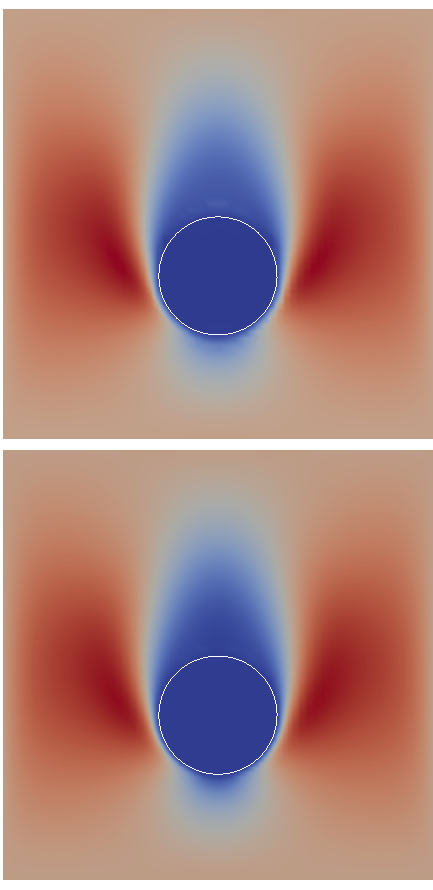}
  \plot{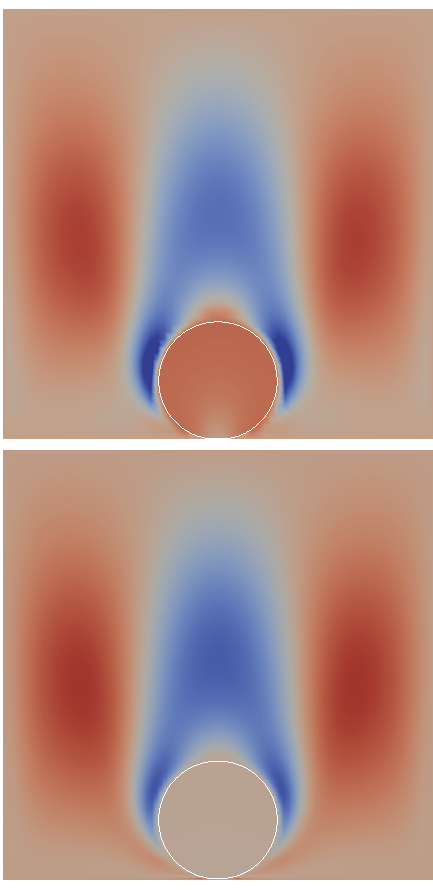}
  \includegraphics[width = \textwidth,
    trim=0cm 10cm 0cm 10cm, clip]{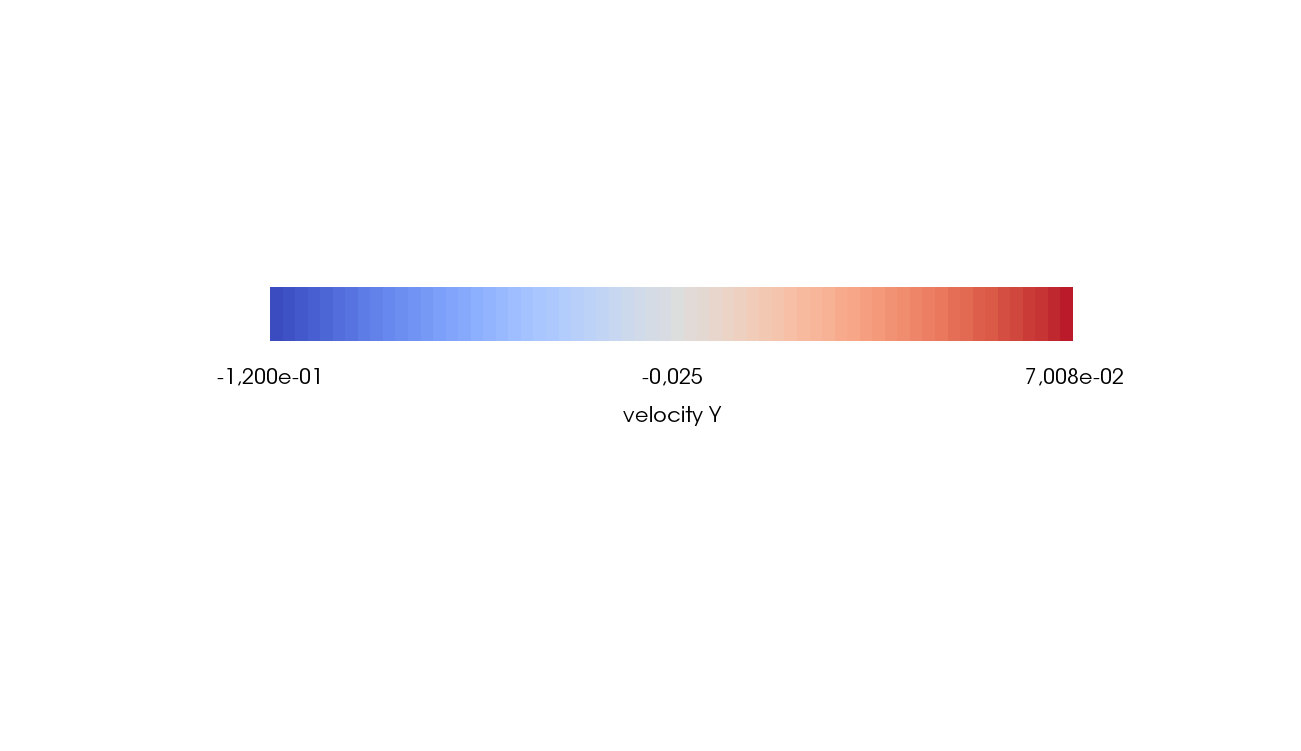}
    \caption{Vertical velocities at time $t=0.29$ (left) and $t=0.489$ (right).
      Top: Uniform mesh. Bottom: Non-uniform mesh.}
\label{fig:plots2}
\end{figure}

Next, we compare the results on mesh level 2 with the locally-refined (non-uniform) mesh described in Section~\ref{sec.geo}.
The vertical velocities of the falling ball and the fluid are illustrated in~\cref{fig:plots2}. The difference between the two results is minor, as long as the ball is still sufficiently far away from the bottom. However, they increase towards contact. This is due to a higher resolution of the pressure peak on the non-uniform mesh, which slows down the fall and
prevents a rebound.

\begin{figure}
  \centering
  % trim={left bottom right top}
  \def\plot#1#2{\includegraphics[width=0.495\textwidth,
    trim=3cm 11cm 3cm 10cm,clip]{pics/sol_p_#1_#2.png}}
  \plot{equi}{478}\hfill\plot{inp}{478}
  \plot{equi}{489}\hfill\plot{inp}{489}
  \plot{equi}{510}\hfill\plot{inp}{510}
  \includegraphics[width = \textwidth,
    trim=0cm 10cm 0cm 10cm, clip]{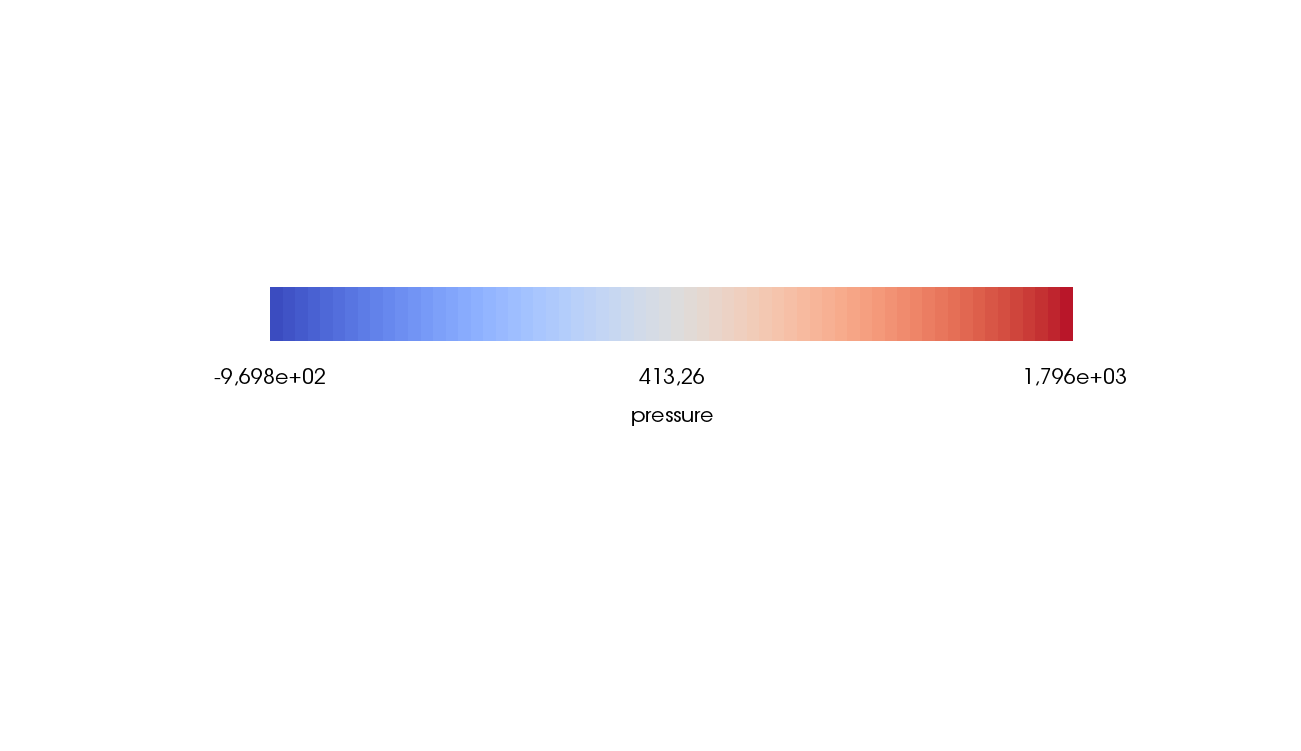}
  \caption{Pressure at time $t=0.478$ (before contact, top), time
      $t=0.489$ (at contact, middle) and $t=0.510$ (after contact, bottom).
      Left: Uniform mesh. Right: Non-uniform mesh. The white horizontal line corresponds to the relaxed distance $\epsilon$. Note that the pressure variable is defined in the fluid domain below and (as an extension) in an additional layer within the solid domain.}
  \label{fig:plots3}
\end{figure}

\begin{figure}
\def\plot#1#2{%
\begin{tikzpicture}
  \small
    \begin{axis}[
        xlabel={Time $t$},
        ylabel={Distance},
        xmin = #1, xmax=#2,
        ymin = 0,
        grid=both,
        width=190pt
    ]
    \addplot[
            thin,
            color=blue
        ]
        table[
            col sep=comma,
            header=false,
            x index=0,
            y index=1
        ] {pics/dist_file1.csv};
        \addlegendentry{Uniform mesh}
    \addplot[
            thin,
            color=red % Farbe der Punkte
        ]
        table[
            col sep=comma,
            header=false,
            x index=0,
            y index=1
        ] {pics/dist_file2.csv};
        \addlegendentry{\small Non-uniform mesh}
    \end{axis}
\end{tikzpicture}
}
\plot{0}{0.6}\hfill\plot{0.47}{0.6}
\caption{Minimal distance between the interface and the bottom boundary against time.}
\label{fig:plots1}
\end{figure}
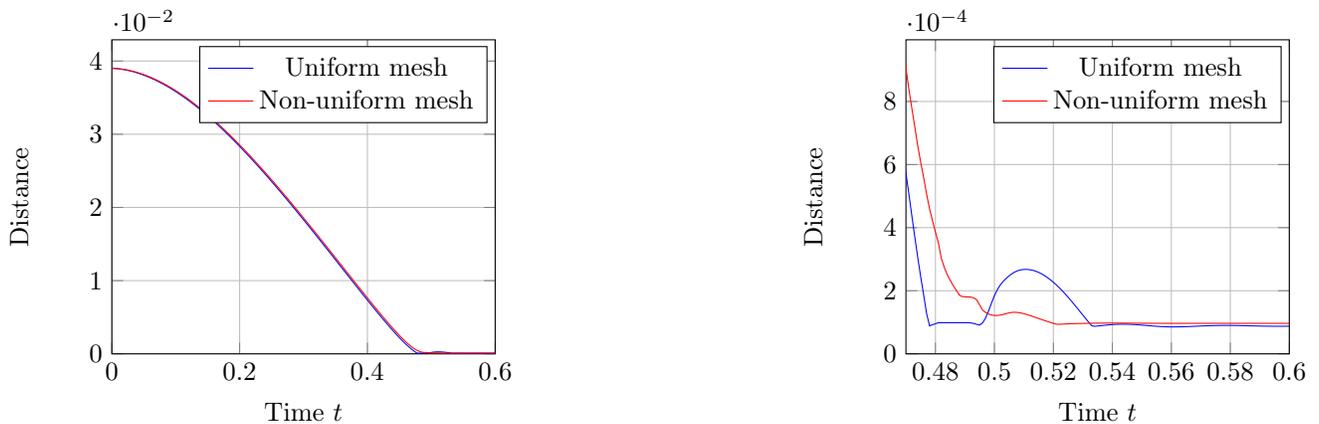

The pressure around the time of contact is displayed in~\cref{fig:plots3}, where the reduction in drop speed for the non-uniform mesh can be observed. Both computations show a stark difference in pressure nearing contact. The distance between the ball and the bottom boundary is presented in~\cref{fig:plots1}. In the zoom at the right part of the figure, we observe that the ball is slowed down strongly before contact on the non-uniform mesh, where the pressure is resolved accurately. In fact, the relaxed contact condition~\eqref{conteps} gets active much later at time $t\approx 0.521$ compared to $t\approx 0.478$ for the uniform mesh. This is expected by arguments related to the Navier-Stokes no-collision paradox~\cite{Hillairet2d, HeslaPhD}.
If the pressure is not resolved sufficiently, as on the uniform mesh, it comes to contact much earlier and a first large rebound as well as a few smaller rebounds are observed.

The further quantities of interest for both cases are given in~\cref{tab:qoi}. We observe that the quantities $t^*, v^*$ and $f^*$, that are defined before the ball is significantly slowed down, are very similar on both meshes. The maximum pressure $p_{bc}$ is approximately by a factor four larger on the non-uniform mesh. On the uniform mesh the maxima of elastic and kinetic energies are similar, which indicates that the kinetic energy is transferred almost completely to elastic energy at the time of impact. On the non-uniform mesh, on the other hand, the elastic energy is very small, as no significant impact takes place. Finally, we observe that the average number of Newton steps required lies below two on both meshes, which shows a very good performance of the nonlinear solver despite the non-smoothness of the equations.

We note that we have performed some additional simulations by varying $\epsilon$ which are, however, not shown in detail here. For larger values of $\epsilon$, the simulations perform well, but the bouncing happens further above the ground (bottom boundary).
If $\epsilon$ is too small, the computation breaks down (as expected) or allows the ball to move past the relaxed contact line.

\begin{table}[tp]
  \label{tab:qoi}
  \caption{Quantities of interest. Top: Uniform mesh.
    Bottom: Non-uniform mesh. All quantities are given in SI units.}
  \centering
  \ifcase0
	\begin{tabular}{*8c}
	\toprule
	DoFs & $t^0$ & $t^*$ & $v^*$ & $f^*$ & $t\tsb{cont}$ & $t\tsb{jump}$ & $h\tsb{jump}$ \\
	\midrule
	37660 & 0.2038 & 0.0613 &$-$0.1034 &$-$4.8917 & 0.2742 & 0.3068 & 0.000267 \\
	55786 & 0.2060 & 0.0614 &$-$0.1031 &$-$4.8698 & 0.3154 & ---    & --- \\
	\bottomrule
\end{tabular}
$ $\\[1ex]
\begin{tabular}{*7c}
	\bottomrule
	DoFs & $\max_t p\tsb{bc}$ & $\max_{Q_f} v_f$ & $\max_t E\tsb{el}$ & $\max_t E_{\text{kin},f}$ & $\max_t E_{\text{kin},s}$ & nNewton \\
	\midrule
	37660 & 575.0429  & 0.1705 & 0.0034 & 0.0090 & 0.0032 & 1.2496 \\
	55786 & 2020.4924 & 0.2417 & 0.0005 & 0.0094 & 0.0031 & 1.3278 \\
	\bottomrule
\end{tabular}
\or
\begin{tabular}{c@{\quad}c@{\quad}c@{\quad}c@{\quad}c@{\quad}c@{\quad}c@{\quad}c}
	\toprule
	DoFs & $t^0$ & $t^*$ & $v^*$ & $f^*$ & $t\tsb{cont}$ & $t\tsb{jump}$ & $h\tsb{jump}$\\
	\midrule
	37660 & 0.2038 & 0.0613 & -0.1034 & -4.8917 & 0.2742 & 0.3070 & 0.000267\\
	55786 & 0.2060 & 0.0614 & -0.1031 & -4.8698 & 0.3154 & --     & -- \\
	\bottomrule\bottomrule
	DoFs & $\max_t p\tsb{bc}$ & $\max_{Q_f} v_f$ & $\max_t E\tsb{el}$ & $\max_t E_{\text{kin},f}$ & $\max_t E_{\text{kin},s}$ & \multicolumn{2}{c}{nNewton}\\
	\midrule
	37660 & 575.0429  & 0.1705 & 0.0034 & 0.0090 & 0.0032 & \multicolumn{2}{c}{1.2516} \\
	55786 & 2020.4924 & 0.2417 & 0.0005 & 0.0094 & 0.0031 & \multicolumn{2}{c}{1.3278} \\
	\bottomrule
\end{tabular}
  \fi
\end{table}

%%%%%%%%%%%%%%%%%%%%%%%%%%%%%%%%%%%%%%%%%%%%%%%%%%%%%%%
\section{Conclusions}
In this work, we proposed a numerical model for formulating fully Eulerian
fluid-structure contact interaction. The setting is described in a monolithic
fashion. As in fully Eulerian systems, the finite element mesh is fixed, the
fluid-structure interface cuts through mesh elements. To treat these
cut elements, a cut finite element approach with ghost penalities is introduced. These parts were worked out in great detail.
In fact, the combination of the fully Eulerian approach for fluid-structure interaction with an unfitted finite element discretization is novel to the knowledge of the authors. A major difference to existing works in a mixed-coordinate framework is the cut finite element formulation for the solid parts in (moving) Eulerian coordinates. This includes, for example, ghost penality terms for the elasticity equations in the solid part, which have to date not been used in the literature within the context of fluid-structure interaction.
Moreover, due to the Eulerian formulation pure hyperbolic terms appear, which are stabilized with SUPG.
Next, the setting is extended to fluid-structure contact interaction.
Here, the no-penetration condition is reformulated as
a non-smooth equality condition, which is inserted into the weak form
of the fluid-structure problem. Finally, the nonlinear solution process is
carried out in a monolithic fashion, which requires directional derivatives
in all three solution variables: velocity, displacements, and pressure. To this end, a semi-smooth line-search Newton scheme was derived with all details
of the Jacobian matrix. In the last section, numerical computations of
an elastic ball in a fluid were conducted. Therein, certain quantities
of interest were evaluated. These showed satisfactory performance in view
of the overall challenging complexity of the problem statement. The results highlight that an accurate resolution of the interface zone, in particular for the pressure variable, is necessary to accurately reproduce the physics of the equations.
Moreover, mesh refinement studies on four meshes were carried out. These show that the coarse mesh is too coarse in terms of the approximation. The three finer meshes show qualitative convergence and thus computational stability of our approach.
In ongoing and future work, computational convergence studies for different time step sizes are being investigated, numerical analysis of a simplified system is underway, other ghost penality implementations are being tested, and finally, iterative linear solvers and preconditioners must be developed.

%%%%%%%%%%%%%%%%%%%%%%%%%%%%%%%%%%%%%%%%%%%%%%%%%%%%%%%%%%%%
\section{Acknowledgements}
Anne-Kathrin Wenske and Marc C.\ Steinbach gratefully acknowledge
  the financial support from the Deutsche Forschungsgemeinschaft
  (DFG, German Research Foundation) -- SFB1463 -- 434502799. Stefan Frei and Thomas Wick
  gratefully acknowledge
  the financial support from the DFG with the grant number 548064929

%%%%%%%%%%%%%%%%%%%%%%%%%%%%%%%%%%%%%%%%%%%%%%%%%%%%%%%%%%%%
%\section{References}

\bibliographystyle{abbrv}
%\bibliography{shorttitles,lit}

\end{document}